\documentclass[leqno]{article}
\usepackage{CJK,CJKnumb,CJKulem,times,dsfont,ifthen,mathrsfs,latexsym,amsfonts,color}
\usepackage{amsmath,amsthm,makeidx,fontenc,amssymb,bm,graphicx,psfrag,listings,curves,extarrows,enumitem}
\usepackage{hyperref}

\allowdisplaybreaks[4] 

\usepackage{geometry}
\usepackage{indentfirst} 
\usepackage{url}
\makeatletter
\def\url@leostyle{%
  \@ifundefined{selectfont}{\def\UrlFont{\sf}}{\def\UrlFont{\small\ttfamily}}}
\makeatother
\usepackage{amssymb}
\makeatletter

\newcommand{\Rmnum}[1]{\expandafter\@slowromancap\romannumeral #1@}
\makeatother

\newtheorem{definition}{Definition}[section]
\newtheorem{theorem}{Theorem}[section]
\newtheorem{lemma}{Lemma}[section]
\newtheorem{corollary}{Corollary}[section]
\newtheorem{proposition}{Proposition}[section]
\newtheorem{remark}{Remark}[section]

\newcommand{\al}{\alpha}
\newcommand{\ga}{\gamma}
\newcommand{\Ga}{\Gamma}
\newcommand{\dl}{\Delta}
\newcommand{\e}{\varepsilon}

\newcommand{\iy}{\infty}
\newcommand{\q}{\theta}

\newcommand{\la}{\lambda}
\newcommand{\vp}{\varphi}
\newcommand{\pa}{\partial}

\newcommand{\ra}{\rightarrow}
\newcommand{\rh}{\rightharpoonup}
\newcommand{\hra}{\hookrightarrow}

\newcommand{\lab}{\label}
\newcommand{\f}{\frac}
\newcommand{\bt}{\begin{theorem}}
\newcommand{\et}{\end{theorem}}
\newcommand{\bl}{\begin{lemma}}
\newcommand{\el}{\end{lemma}}
\newcommand{\bd}{\begin{definition}}
\newcommand{\ed}{\end{definition}}
\newcommand{\bc}{\begin{corollary}}
\newcommand{\ec}{\end{corollary}}
\newcommand{\bp}{\begin{proof}}
\newcommand{\ep}{\end{proof}}
\newcommand{\bx}{\begin{example}}
\newcommand{\ex}{\end{example}}
\newcommand{\bi}{\begin{exercise}}
\newcommand{\ei}{\end{exercise}}
\newcommand{\br}{\begin{remark}}
\newcommand{\er}{\end{remark}}
\newcommand{\be}{\begin{equation}}
\newcommand{\ee}{\end{equation}}
\newcommand{\bal}{\begin{align}}
\newcommand{\bn}{\begin{enumerate}}
\newcommand{\en}{\end{enumerate}}
\newcommand{\ba}{\begin{align}}
\newcommand{\ea}{\begin{align}}
\newcommand{\bg}{\begin{align*}}
\newcommand{\eg}{\end{align*}}
\newcommand{\bcs}{\begin{cases}}
\newcommand{\ecs}{\end{cases}}
\newcommand{\AR}{{\cal A}}

\newcommand{\CR}{{\cal C}}
\newcommand{\DR}{{\cal D}}
\newcommand{\ER}{{\cal E}}
\newcommand{\FR}{{\cal F}}
\newcommand{\GR}{{\cal G}}
\newcommand{\HR}{{\cal H}}

\newcommand{\OR}{{\cal O}}
\newcommand{\QR}{{\cal Q}} 
\newcommand{\SR}{{\cal S}} 

\newcommand{\XR}{{\cal X}}

\newcommand{\N}{{\mathbb N}}

\newcommand{\R}{{\mathbb R}}
\newcommand{\RN}{{\mathbb R^N}}  
\newcommand{\bean}{\begin{eqnarray*}}
\newcommand{\eean}{\end{eqnarray*}}

\newcommand{\sbr}[1]{\left(#1\right)}
\newcommand{\mbr}[1]{\left[#1\right]}
\newcommand{\lbr}[1]{\left\{#1\right\}}
\newcommand{\rd}{\mathrm d}

\newcommand{\ind}{\mathrm{Ind}}
\newcommand{\s}{\star}

\newcommand{\nm}[1]{\Vert #1 \Vert}
\newcommand{\nt}[1]{\|#1\|_2}
\newcommand{\gu}[1]{|\nabla u|^}
\newcommand{\n}{{\mu_n}}
\newcommand{\nj}{{\mu_n^j}}

\setitemize{itemindent=38pt,leftmargin=0pt,itemsep=-0.4ex,listparindent=26pt,partopsep=0pt,parsep=0.5ex,topsep=-0.25ex}

\numberwithin{equation}{section}

\begin{document}
\theoremstyle{plain}

\title{\bf Quasilinear Schr\"odinger equations: ground state and infinitely many normalized solutions
	\thanks{This work is supported by NSFC(11771234,12026227);  E-mails: li-hw17@mails.tsinghua.edu.cn\quad\& \quad zou-wm@mail.tsinghua.edu.cn}
	}

\date{}
\author{
{\bf Houwang Li$^1$\;\&\;Wenming Zou$^2$}\\
\footnotesize \it 1. Department of Mathematical Sciences, Tsinghua University, Beijing 100084, China.\\
\footnotesize \it 2. Department of Mathematical Sciences, Tsinghua University, Beijing 100084, China.}

\maketitle
\begin{center}
\begin{minipage}{120mm}
\begin{center}{\bf Abstract }\end{center}
In the present paper, we study the normalized solutions for the following quasilinear Schr\"odinger equations:
	$$-\dl u-u\dl u^2+\la u=|u|^{p-2}u \quad \text{in}~\RN,$$
with prescribed mass
	$$\int_{\RN} u^2=a^2.$$
We first consider the mass-supercritical case $p>4+\f{4}{N}$, which has not been studied before. 
By using a perturbation method, we succeed to prove the existence of ground state normalized solutions, and by applying the  index theory, we 
obtain the existence of infinitely many normalized solutions. Then we turn to study the mass-critical case, i.e., $p=4+\f{4}{N}$, and obtain 
some new existence results. Moreover, we also  observe  a concentration behavior of the ground state solutions.
\vskip0.23in



\vskip0.23in

\end{minipage}
\end{center}

\vskip0.23in
\section{Introduction}
We consider the equation
\be\lab{mainequ}
	-\dl u-u\dl u^2+\la u=|u|^{p-2}u, \quad \text{in}~\RN,
\ee
which is usually called Modified Nonlinear Schr\"odinger equation. Such type of equations appear 
as a standing wave version of the following Schr\"odinger equations,
\be\lab{parabolicequ}
	\left\{ \begin{aligned}
	& i\partial_t \phi+\dl \phi+\phi\dl(|\phi|^2)+|\phi|^{p-2}\phi=0,\quad \text{in}~\R^+\times\RN,\\
	& \phi(0,x)=\phi_0(x),\quad \text{in}~\RN.
	\end{aligned} \right.
\ee
It is well known that the above Schr\"odinger equations model many phenomena in mathematical physics, for instance in the theory of Heisenberg 
ferromagnets and magnons \cite{bg-1,bg-2,bg-3}, in models of superfluid films in fluid mechanics and plasma physics \cite{bg-4,bg-5,bg-6}, 
in dissipative quantum mechanics \cite{bg-7}, and in condensed matter theory \cite{bg-8}, which have received considerable attention 
in mathematical analysis during the last two decades.

\vskip0.1in

	In recent years, the search  for  the solution with prescribed mass has became a hot direction, that is to find $u$ such that
\be\lab{massequ}
	\left\{ \begin{aligned}
	&-\dl u-u\dl u^2+\la u=|u|^{p-2}u, \quad \text{in}~\RN,\\
	&\int_\RN |u|^2\rd x=a,
	\end{aligned}\right.
\ee
with $\la$ appearing as Lagrange multiplier. From the view of physics, prescribed mass represents the law of conservation of mass, 
so it seems to be great meaningful to study such solutions. Solutions of prescribed mass are often referred to as normalized solutions, and
the present paper is devoted to such solutions i.e., the solution $u$ of \eqref{massequ} with a Lagrange multiplier $\la\in\R$.

\vskip0.1in

	The existence of normalized solutions to the semilinear Schr\"odinger equation
\be\lab{semilinearequ}
	-\dl u+\la u=g(u), \quad \text{in}~\RN,
\ee
has been widely studied recently. Mathematically, to obtain the normalized solutions, one needs  to consider
the corresponding energy functional on a $L^2$ sphere, which has particular difficulties: the weak limit of the Palais-Smale sequence 
may be not contained in the $L^2$ sphere (even in the radial case), and the Palais-Smale sequence does not even need to be bounded.
So the study of normalized solutions of \eqref{semilinearequ} is much more complicated than 
the study of \eqref{semilinearequ} with prescribed $\la\in\R$. Fortunately, in \cite{Jeanjean=1997}, using an  auxiliary 
functional and a mini-max theorem from \cite{Ghoussoub=1993}, L. Jeanjean obtained a normalized solution of \eqref{semilinearequ}. 
The existence of infinitely many normalized solutions of \eqref{semilinearequ} was later proved by T. Bartsch and S. de Valeriola in 
\cite{Bartsch-Valeriola=2013} using a new linking geometry for the auxiliary functional (see also the papers by T. Bartsch and N. Soave 
\cite{BartschSoave=2017}). After that, N. Ikoma and K. Tanaka \cite{Ikoma=AdvDE=2019} constructed a deformation theorem suitable for
the auxiliary functional, and then obtained infinitely many normalized solutions of \eqref{semilinearequ} through Krasnoselskii index under
a weaker condition on $g(u)$. Soon later, L. Jeanjean and S. S. Lu \cite{Jeanjean-Lu=CVPDE=2020} obtained infinitely many normalized 
solutions of \eqref{semilinearequ} under  a totally different assumption on $g(u)$ which permits $g(u)$ to be just continuous. As for the
least energy normalized solutions, N. Soave in \cite{Soave=JDE=2020,Soave=JFA=2020}, by restraining the energy functional on 
a smaller manifold, obtained the existence of ground state normalized solutions with $g(u)=|u|^{p-2}u+\mu|u|^{q-2}u$. For more 
results on normalized solutions for scalar equations and systems, we refer to \cite{Zhong-Zou-Bartsch=2020,BS=CVPDE=2019,BJS=JMPA=2016,
BartschJeanjean=2018,GouJeanjean=2018,GouJeanjean=2016,Li-Zou,Bellazzini-Jeanjean-Luo=2013}.

\vskip0.1in
	Now we come back to the Modified Nonlinear Schr\"odinger equation \eqref{mainequ}. When considering \eqref{mainequ} with $\la\in\R$ fixed, 
one always study the functional
\be
	E_\la(u):=\f{1}{2}\int_\RN(\gu^2+\la|u|^2) +\int_\RN |u|^2\gu^2-\f{1}{p}\int_\RN|u|^p,
\ee
on the space 
	$$\HR=\lbr{u\in W^{1,2}(\RN):\int_\RN |u|^2\gu^2<+\iy}.$$
It is easy to check that $u$ is a weak solution of \eqref{mainequ} if and only if
	$$E_\la'(u)\phi=\lim_{t\to0^+}\f{E_\la(u+t\phi)-E_\la(u)}{t}=0,$$
for every $\phi\in \CR_0^\iy(\RN)$. We recall, see \cite{Liu-Wang-Wang=JDE=2003} for example, 
that the value 
	$$22^*=\begin{cases}&\f{4N}{N-2},\quad N\ge3,\\ &+\iy,\quad N\le2 \end{cases}$$ 
corresponds to a critical exponent. Compared to equation \eqref{semilinearequ}, the search 
of solutions of \eqref{mainequ} presents a major difficulty: the functional associated with the term $u\dl u^2$
	$$V(u)=\int_\RN |u|^2\gu^2$$
is non-differentiable in $\HR$ when $N\ge2$. To overcome this difficulty, various arguments have been developed, such as 
the minimizition methods \cite{Poppenberg-Schmitt-Wang=CVPDE=2002} where the non-differentiability of $E_\la$ does 
not come into play, the methods of a Nehari manifold approach \cite{Liu-Liu-Wang=CVPDE=2013}, the methods of changing variables 
\cite{Liu-Wang-Wang=JDE=2003,Colin-Jeanjean=NA=2004} which transform problem \eqref{mainequ} into a semilinear one \eqref{semilinearequ},
and a perturbation method in a series of papers \cite{Liu-Liu-Wang=ProAMS=2013,Liu-Wang=JDE=2014,Liu-Liu-Wang=JDE=2013} 
which recovers the differentiability by considering a perturbed functional on a smaller function space.

\vskip0.1in
	However, when considering the normalized solution problem \eqref{massequ}, one would find that the methods of Nehari manifold approach 
and changing variables are no longer applicable, since the parameter $\la$ is unknown and the $L^2$-norm $\nm{u}_2$ must 
be equal to a given number. So there are very few results on problem \eqref{massequ}. Formally, a normalized solution of \eqref{massequ} 
can be obtained as a critical point of 
\be\lab{energy}
	I(u):=\f{1}{2}\int_\RN \gu^2 +\int_\RN |u|^2\gu^2-\f{1}{p}\int_\RN|u|^p
\ee
on the set
\be\lab{massball}
	\tilde\SR(a):=\lbr{u\in \HR: \int_\RN|u|^2=a},
\ee
that is, a normalized solution of \eqref{massequ} is a $u\in \tilde\SR(a)$ such that there exists a $\la\in\R$ satisfing
\be
	\int_\RN \nabla u\cdot\nabla\phi +2\int_\RN (u\phi\gu^2+|u|^2\nabla u\cdot\nabla\phi)+\la\int_\RN u\phi-\int_\RN|u|^{p-2}u\phi=0,
\ee
for any $\phi\in\CR_0^\iy(\RN)$. To proceed our paper, we introduce a sharp Gagliardo-Nirenberg inequality \cite{Agueh=2008}:
\be\lab{SGNine}
	\int_\RN |u|^{\f{p}{2}}\le \f{C(p,N)}{\nm{Q_p}_1^\f{p-2}{N+2}}\sbr{\int_\RN|u|}^{\f{4N-(N-2)p}{2(N+2)}}
		\sbr{\int_\RN|\nabla u|^2}^{\f{N(p-2)}{2(N+2)}},\quad \forall~u\in \ER^1,
\ee
where $2<p<22^*$,
	$$C(p,N)=\f{p(N+2) }{\mbr{4N-(N-2)p}^{\f{4-N(p-2)}{2(N+2)}}\mbr{2N(p-2)}^{\f{N(p-2)}{2(N+2)}} },$$
and 
	$$\ER^q:=\lbr{u\in L^q(\RN):\nabla u\in L^2(\RN)},$$
with norm $\nm{u}_{\ER^q}:=\nt{\nabla u}+\nm{u}_q$. It is well known that $\ER^q$ is a reflexive Banach space when $1<q<\iy$, 
and for Embeeding theorems and more related properties we refer to \cite{Kuzin-Pohozaev=1997}. 
Moreover, $Q_p$ optimizes \eqref{SGNine} and is the unique nonnegative radially symmetric solution of the following 
equation \cite{Serrin-Tang=2000}:
\be\lab{GNequ1}
	-\Delta u+1=u^{\f{p}{2}-1},\quad \text{in}~\RN.
\ee
Strictly speaking, it has been proved in \cite[Theorem 1.3]{Serrin-Tang=2000} that $Q_p$ has a compact support in $\RN$ 
	and exactly satisfies a Dirichlet-Neumann free boundary problem. Namely, there exists an  $R>0$ such that $Q_p$ is the unique
	positive solution of 
	\be\lab{GNequ2}
		\left\{ \begin{aligned}
		&-\Delta u+1=u^{\f{p}{2}-1},\quad \text{in}~B_R,\\
		&u=\f{\pa u}{\pa n}=0,\quad\text{on}~\pa B_R.
		\end{aligned}\right.
	\ee
	In what follows, if we say that $u$ is a nonnegative solution of \eqref{GNequ1}, then we  mean that $u$ is a 
	solution of \eqref{GNequ2}.
By replacing $u$ with $u^2$ in \eqref{SGNine}, one immediately obtain the following Gagliardo-Nirenberg-type inequality,
\be\lab{GNine}
	\int_\RN|u|^p\le \f{C(p,N)}{\nm{Q_p}_1^\f{p-2}{N+2}}\sbr{\int_\RN|u|^2}^{\f{4N-p(N-2)}{2(N+2)}}
			\sbr{4\int_\RN |u|^2\gu{u}2}^{\f{N(p-2)}{2(N+2)}}.
\ee

	Now we collect some known results about normalized solutions of \eqref{massequ}. First, to avoid 
the nondifferentiability of $V(u)$, 
M. Colin, L. Jeanjean, M. Squassina \cite{Colin-Jeanjean-Squassina=2010} and L. Jeanjean, T. J. Luo \cite{Jeanjean-Luo=2013} 
considered the minimizition problem
	$$\tilde m(a)=\inf_{u\in \tilde\SR(a)} I(u)$$
with $2<p\le4+\f{4}{N}$. Using ineqality \eqref{GNine}, one can find that $\tilde m(a)>-\iy$ when $2<p<4+\f{4}{N}$ and $\tilde m(a)=-\iy$ 
when $p>4+\f{4}{N}$, since
	$$\f{N(p-2)}{2(N+2)}<1\quad \text{ if and only if } \quad p<4+\f{4}{N}.$$
These considerations show that the exponent $4+\f{4}{N}$  for   \eqref{massequ} plays  the role of $2+\f{4}{N}$ in \eqref{semilinearequ}. 
After that, X. Y. Zeng and Y. M. Zhang \cite{Zeng-Zhang=2018} studied the existence and asymptotic behavior of the minimiziers to  
	$$\inf_{u\in\tilde\SR(a)}I(u)+\int_\RN a(x)|u|^2,$$
where $a(x)$ is an infinite pontential well.  In addition to these minimizition approaches, L. Jeanjean, T. J. Luo and Z. Q. Wang  
\cite{Jeanjean-Luo-Wang=JDE=2015} obtained another mountain-pass type normalized solution of \eqref{massequ} through the perturbation method.
We remark that all of these results on normalized solution of \eqref{massequ} only considered the mass-subcritical or 
mass-critical case, i.e.,  $2<p\le4+\f{4}{N}$.

\vskip0.2in

In this paper, we consider the mass-critical and mass-supercritical case, i.e.,  $p\ge4+\f{4}{N}$.  To the best of our knowledge, the case of  mass-supercritical has not been 
considered before.  Actually, we obtain 
\bt\lab{thm1}
	Assume that one of the following conditons holds:
	\begin{itemize}[fullwidth,itemindent=2em]
	\item[(H1)]	$N=1,2$, $p>4+\f{4}{N}$, $a>0$;
	\item[(H2)]	$N=3$, $4+\f{4}{N}<p<2^*$, $a>0$.
	\end{itemize}
	Then there exists a radially symmetric positive ground state normalized solution
		$u\in W^{1,2}(\RN)\cap L^\iy(\RN)$ of \eqref{massequ} in the sense that
	$$I(u)=\inf\lbr{I(v): v\in \tilde\SR(a), I|_{\tilde\SR(a)}'(v)=0, v\neq0}.$$
\et

\bt\lab{thm2}
	Assume that one of the following conditons holds
	\begin{itemize}[fullwidth,itemindent=2em]
	\item[(H1)']	$N=2$, $p>4+\f{4}{N}$, $a>0$,
	\item[(H2)]		$N=3$, $4+\f{4}{N}<p<2^*$, $a>0$.
	\end{itemize}
	Then there exists a sequence of normalized solutions 
	$u^j\in W^{1,2}(\RN)\cap L^\iy(\RN)$ of \eqref{massequ} with increasing energy $I(u^j)\to+\iy$.
\et
\br
	\begin{itemize}[fullwidth,itemindent=0em]
	\item[(1)]	We state that the dimension is limited due to a lemma limitition used to control the Lagrange multipliers, 
	see Lemma \ref{Liouville1} and Remark \ref{remark2}.
	\item[(2)]	The difference between Theorem \ref{thm1} and Theorem \ref{thm2} is that we can not prove the existence of 
	infinitely many solutions when $N=1$, because the failure of the compact embedding $W^{1,2}(\R)\hra\hra L^q(\R)$ for 
	$2<q<2^*$. However when considering the ground state, we are able to recover the compactness of bounded sequences using 
	the symmetric decreasing arrangement, due
	to the advantage of the associated minimizition $m_\mu(a)$ defined in \eqref{min}.
	\end{itemize}
\er

	Now we turn to the mass-critical case, i.e.,  $p=4+\f{4}{N}$. Let $a_*=\nm{Q_{4+\f{4}{N}}}_1$.
\bt\lab{thm3}
	Assume that one of the following conditons holds:
	\begin{itemize}[fullwidth,itemindent=2em]
	\item[(H3)]	$N\le3$, $p=4+\f{4}{N}$, $a>a_*$; 
	\item[(H4)]	$N\ge4$, $p=4+\f{4}{N}$, $a_*<a<\sbr{\f{N-2}{N-2-\f{4}{N}}}^{\f{N}{2}}a_*$, 
	\end{itemize}
	Then there exists a radially symmetric positive ground state normalized solution
		$u\in W^{1,2}(\RN)\cap L^\iy(\RN)$ of \eqref{massequ} in the sense that
	$$I(u)=\inf\lbr{I(v): v\in \tilde\SR(a), I|_{\tilde\SR(a)}'(v)=0, v\neq0}.$$
\et
\br
	In a very recent paper \cite{Ye=JMAA=2021}, H. Y. Ye and Y. Y. Yu obtained the existence of ground state normalized solution 
	of \eqref{massequ} under assumption $(H3)$. As one can see, although Theorem \ref{thm3} contains their existence result, 
	the method we used in the current paper  is totally different from theirs, while as they said in \cite[Remark 1.3]{Ye=JMAA=2021}, they are unable to 
	handle the case $N\ge4$. Moreover, they also consider a asymptotic behavior, but our Theorem \ref{thm4} is more 
	accurate, since we give a discription of $u_n$ when $a\to a_*$.
\er

	We observe that when $p=4+\f{4}{N}$, the value $a_*$ is a threshold of the existence of normalized solution of \eqref{massequ}.
Actually, we have
\begin{proposition}\lab{prop1}
	Let $p=4+\f{4}{N}$ and $N\ge1$. Then
	\begin{itemize}[fullwidth,itemindent=2em]
	\item[(1)]	$\tilde m(a)=\begin{cases} &0,\quad 0< a\le a_*,\\
						&-\iy,\quad a>a_*.\end{cases}$
	\item[(2)]	\eqref{massequ} has no solutions for any $0<a\le a_*$.
	\item[(3)]	\eqref{massequ} has at least one radially symmetric positive solution for $a>a_*$ and $a$ is close to $a_*$.
	\end{itemize}
\end{proposition}
\br\lab{remark1}
	We state that (1) is a direct conclusion of \cite[Theorem 1.9]{Colin-Jeanjean-Squassina=2010}, and (3) is a direct conclusion
	of Theorem \ref{thm3} above. Now we prove (2). Since $u$ is a solution of \eqref{massequ}, there holds (see Lemma \ref{lempho})
		$$\int_\RN\gu{u}2+(2+N)\int_\RN|u|^2\gu{u}2-\f{N(2+N)}{4(N+1)}\int_\RN|u|^{4+\f{4}{N}}=0.$$
	Combining with \eqref{GNine}, we obtain
		$$\int_\RN\gu{u}2+(2+N)\int_\RN|u|^2\gu{u}2\le (2+N)\sbr{\f{a}{a_*}}^{\f{2}{N}}\int_\RN|u|^2\gu{u}2,$$
	from which we get $u=0$ for any $0<a\le a_*$, a contradiction since $\nt{u}=a$.
\er

	Inspired by Proposition \ref{prop1}, we enlighten a concentration behavior of the radially symmetric positive solution of 
\eqref{massequ} when $p=4+\f{4}{N}$ and $a\to a_*$.

\bt\lab{thm4}
	Let $p=4+\f{4}{N}$, $N\ge1$, and let $u_n$ be a radially symmetric positive solution of \eqref{massequ} for $a=a_n$ with 
	$a_n>a_*$ and $a_n\to a_*$. Then there exists a sequence $y_n\in\RN$ such that up to a subsequence,
	\be\lab{concentration}
		\mbr{\sbr{\f{Na_*}{N}}^{\f{1}{2+N}}\e_n}^{N} u_n^2\left(\sbr{\f{Na_*}{N}}^{\f{1}{2+N}}\e_nx +\e_n y_n\right)
		\to Q_{4+\f{4}{N}}\quad \text{in}~L^q(\RN)
	\ee
	for $1\le q< 2^*$, where 
		$$\e_n=\sbr{\int_\RN u_n^2|\nabla u_n|^2}^{-(2+N)}\to0.$$
\et

\br
	Theorem \ref{thm4} gives a description of radially symmetric positive solution of \eqref{massequ} as the mass $a_n$ approaches 
	to $a_*$ from above. Roughly speaking, it shows that for $N$ large enough, we have
	$$u_n(x)=\mbr{\sbr{\f{Na_*}{N}}^{\f{1}{2+N}}\e_n}^{-\f{N}{2}} 
		Q_{4+\f{4}{N}}\left(\sbr{\f{Na_*}{N}}^{-\f{1}{2+N}}\e_n^{-1}(x-\e_n^{-1}y_n) \right).$$
\er

	The paper is organized as follows. In Section 2, we give perturbation settings and an important lemma.
In section 3.1, we give some properties of the associated Pohozaev manifold. In section 3.2 and 3.3, we prove the 
existence of ground state and infinitely many critical points for perturbed funtional. In section 4, we study the convergence 
of the critical points the perturbated funtional as $\mu\to0^+$. And the Theorem \ref{thm1} for $N=1$ is proved in section 
\ref{groundstatesolution}; the Theorem \ref{thm1} for $N\ge2$ and Theorem \ref{thm2} are proved in section 4. 
Finally, in section 5, we study the mass-critical case, and prove Theorems \ref{thm3}, \ref{thm4}. In the Appendix, we prove 
some valuable results.

\vskip0.1in

	Throughtout the paper, we use standard notations. For simplicity, we write $\int_\RN f$ to mean the Lebesgue integral of 
$f(x)$ over $\RN$. $\nm{\cdot}_p$ denotes the standard norm of $L^p(\RN)$;
We use ``$\to$'' and ``$\rh$'' to denote the strong and weak convergence in the related function space respectively; 
$C,C_1,C_2,\cdots$ will denote positive constants unless specified.

\vskip0.23in
\section{Preliminary}

\subsection{Perturbation setting}\lab{sec1}
	Let $I(u)$ be defined by \eqref{energy}. Observe that when $N=1$, $I(u)$ is of calss $\CR^1$ in $W^{1,2}(\R)$, 
so there is no need to perturb $I(u)$, and in this case 
the proof will be stated separately in the last of Section \ref{groundstatesolution}.
Thus we assume $N\ge2$. To avoid the non-differentiability, we define for $\mu\in(0,1]$,
\be\lab{energy1}
	I_\mu(u):=\f{\mu}{\q}\int_\RN|\nabla u|^\q +I(u)
\ee
on the space $\XR:=W^{1,\q}(\RN)\cap W^{1,2}(\RN)$ for some fixed $\q$ satisfying
	$$\f{4N}{N+2}<\q<\min\lbr{\f{4N+4}{N+2},N}\quad \text{when}~N\ge3$$
and
	$$2<\q<3\quad \text{when}~N=2.$$
Then $\XR$ is a reflexive Banach space. And Lemma \ref{A1} implies $I_\mu\in\CR^1(\XR)$. 
We will consider $I_\mu$ on the constraint
\be
	\SR(a):=\lbr{u\in\XR:\int_\RN|u|^2=a}.
\ee
Recalling the $L^2$-norm preserved transform \cite{Jeanjean=1997}
	$$u\in\SR(a)\mapsto s\s u(x)=e^{\f{N}{2}s}u(e^s x)\in\SR(a),$$
we define 
\be\lab{phoequ}
	\begin{aligned}
	&Q_\mu(u)\\
	&:=\f{\rd}{\rd s}\big|_{s=0} I_\mu(s\s u)\\
			&=(1+\ga_\q)\mu\int_\RN\gu{u}\q+\int_\RN\gu{u}2+(2+N)\int_\RN|u|^2\gu{u}2-\ga_p\int_\RN|u|^p,
	\end{aligned}
\ee
where $\ga_p=\f{N(p-2)}{2p}$. And again Lemma \ref{A1} implies $Q_\mu\in\CR^1(\XR)$. Then we define the manifold
\be
	\QR_\mu(a):=\lbr{u\in\SR(a):Q_\mu(u)=0}.
\ee
We observed that
\bl\lab{lempho}
	Any critical point $u$ of $I_\mu|_{\SR(a)}$ is contained in $\QR_\mu(a)$.
\el
\bp
By \cite[Lemma 3]{Berestycki-Lions=1983-2}, there exists a $\la\in\R$ such that
\be\lab{tem1}
	I_\mu'(u)+\la u=0\quad \text{in}~\XR^*.
\ee
On the one hand, testing \eqref{tem1} with $x\cdot\nabla u$, see \cite[Proposition 1]{Berestycki-Lions=1983-1} for details, we obtain
\be\lab{tem2}
	\begin{aligned}
	0&=\f{\q-N}{\q}\mu\int_\RN\gu{u}\q+\f{2-N}{2}\int_\RN\gu{u}2+(2-N)\int_\RN|u|^2\gu{u}2\\
		&\quad~~ +\f{N}{p}\int_\RN|u|^p-\f{N}{2}\la\int_\RN|u|^2.
	\end{aligned}
\ee
On the other hand, testing \eqref{tem1} with $u$, we obtain
\be\lab{tem3}
	0=\mu\int_\RN\gu{u}\q+\int_\RN\gu{u}2+4\int_\RN|u|^2\gu{u}2-\int_\RN|u|^p+\la\int_\RN|u|^2.
\ee
Combining \eqref{tem2} and \eqref{tem3}, we have $Q_\mu(u)=0$. Then $u\in\QR_\mu(a)$.
\ep

\subsection{An important lemma}

	We need the following result, which are crucially used to control the possible values of  the Lagrange parameters.
\bl\lab{Liouville1}
	Suppose $u\neq0$ is a critical point of $I_\mu|_{\SR(a)}$ with $0\le\mu\le1$, that is there exists a $\la\in\R$ such that 
		$$I_\mu'(u)+\la u=0\quad \text{in}~\XR^*.$$
	And assume that one of the following conditions holds
	\begin{itemize}[fullwidth,itemindent=2em]
	\item[(a)]	$1\le N\le3$, $4+\f{4}{N}\le p\le 2^*$, $a>0$,
	\item[(b)]	$N\ge4$, $p=4+\f{4}{N}$, $0<a<\sbr{\f{N-2}{N-2-\f{4}{N}}}^{\f{N}{2}}a_*$,
	\end{itemize}
	then $\la>0$.
\el
\bp
	By combining $Q_\mu(u)=0$ and \eqref{thm3}, we obtain
		$$\begin{aligned}
			\f{\la N(p-2)}{2p}a &=\sbr{1+\f{N(p-\q)}{p\q}}\mu\int_\RN|\nabla u|^\q+\f{2N-(N-2)p}{2p}\int_\RN|\nabla u|^2\\
			 &\quad~ +\f{4N-(N-2)p}{2p}\int_\RN u^2|\nabla u|^2.
		\end{aligned}$$
So if condition (a) holds, we immediately get $\la>0$. Now suppose condition (b) holds. 
Again from $Q_\mu(u)=0$ and \eqref{tem3}, and using inequality \eqref{GNine}, we obtain
	$$\begin{aligned}
	\la a &=\f{N(\q-2)}{2\q}\mu\int_\RN|\nabla u|^\q+(N-2)\int_\RN u^2|\nabla u|^2
			-\f{N^2-2N-4}{4(N+1)}\int_\RN|u|^{4+\f{4}{N}}\\
		 &\ge \mbr{(N-2)-(N-2-\f{4}{N})\sbr{\f{a}{a_*}}^{\f{2}{N}} }\int_\RN u^2|\nabla u|^2\\
		 &>0,
	\end{aligned}$$
which gives $\la>0$.
\ep

\vskip0.23in
\section{The critical points of perturbed functional}\lab{sec3}
	
	In this whole section, we assume $p>4+\f{4}{N}$.

\subsection{Properties of \texorpdfstring{$\QR_\mu(a)$}{} }

\bl
	Let $0<\mu\le1$, then $\QR_\mu(a)$ is a $\CR^1$-submanifold of codimension 1 in $\SR(a)$, hence a $\CR^1$-submanifold of 
codimension 2 in $\XR$.
\el
\bp
	As a subset of $\XR$, the set $\QR_\mu(a)$ is defined by the two equations $G(u)=0$, $Q_\mu(u)=0$, where
	$$G(u)=a-\int_\RN |u|^2,$$
and clearly $G\in\CR^1(\XR)$. We have to check that
\be
	\rd(Q_\mu,G):\XR\to\R^2\quad \text{is surjective}.
\ee
If this is not ture, $\rd Q_\mu(u)$ and $\rd G(u)$ are linearly dependent, i.e., there exists $\nu\in\R$ such that
\be
	\begin{aligned}
	& \quad~~\q(1+\ga_\q)\mu\int_\RN\gu{u}{\q-2}\nabla u\cdot\nabla\phi +2\int_\RN\nabla u\cdot\nabla \phi \\
	& +(2+N)2\int_\RN(|u|^2\nabla u\cdot\nabla \phi+u\phi\gu{u}2)-p\ga_p\int_\RN|u|^{p-2}u\phi=2\nu\int_\RN u\phi,
	\end{aligned}
\ee
for any $\phi\in\XR$.
Similar as Lemma \ref{lempho}, taking $\phi=x\cdot\nabla u$ and $\phi=u$, we obtain
\be
	\q(1+\ga_\q)^2\mu\int_\RN\gu{u}\q+2\int_\RN\gu{u}2+(2+N)^2\int_\RN|u|^2\gu{u}2-p\ga_p^2\int_\RN|u|^p=0.
\ee
Since $Q_\mu(u)=0$, we get
\be
	\begin{aligned}
	& (p\ga_p-\q-\q\ga_\q)(1+\ga_\q)\mu\int_\RN\gu{u}\q+(p\ga_p-2)\int_\RN\gu{u}2\\
	&\quad~~ +(p\ga_p-2-N)(2+N)\int_\RN|u|^2\gu{u}2=0,
	\end{aligned}
\ee
which means $u=0$ since $p\ga_p>\q+\q\ga_\q$ and $p\ga_p>2+N$. That contradics with $u\in\SR(a)$.
\ep

Now, we prove the following important lemma.
\bl\lab{structure}
	For any $0<\mu\le1$ and any $u\in\XR\setminus\{0\}$, the following statements hold.
	\begin{itemize}[fullwidth,itemindent=0em]
	\item[(1)]	There exists a unique number $s_\mu(u)\in\R$ such that $Q_\mu(s_\mu(u)\s u)=0$.
	\item[(2)]	$I_\mu(s\s u)$ is strictly increasing in $s\in(-\iy,s_\mu(u))$ and is strictly decreasing 
	in $s\in(s_\mu(u),+\iy)$, and
		 $$\lim_{s\to-\iy}I_\mu(s\s u)=0^+,\quad \lim_{s\to+\iy}I_\mu(s\s u)=-\iy,\quad I_\mu(s_\mu(u)\s u)>0.$$
	\item[(3)]	$s_\mu(u)<0$ if and only if $Q_\mu(u)<0$.
	\item[(4)]  The map $u\in\XR\setminus\{0\} \mapsto s_\mu(u)\in\R$ is of class $\CR^1$.
	\item[(5)]  $s_\mu(u)$ is an even function with respect to $u\in \XR\setminus\{0\}$.
	\end{itemize}
\el
\bp
\begin{itemize}[fullwidth,itemindent=0em]
\item[(1)] By direct computation, one can check that 
\be
	\begin{aligned}
	Q_\mu(s\s u):&=\f{\rd}{\rd s} I_\mu(s\s u)\\
			&=(1+\ga_\q)\mu e^{\q(1+\ga_\q)s}\int_\RN\gu{u}\q+e^{2s}\int_\RN\gu{u}2\\
			&\quad~~ +(2+N)e^{(2+N)s}\int_\RN|u|^2\gu{u}2-\ga_pe^{p\ga_ps}\int_\RN|u|^p\\
			&=e^{p\ga_ps}\big[ (1+\ga_\q)\mu e^{-(p\ga_p-\q-\q\ga_\q)s}\int_\RN\gu{u}\q+e^{-(p\ga_p-2)s}\int_\RN\gu{u}2\\
			&\quad~~ +(2+N)e^{-(p\ga_p-2-N)s}\int_\RN|u|^2\gu{u}2-\ga_p\int_\RN|u|^p \big]
	\end{aligned}
\ee
Since $p\ga_p>\q+\q\ga_\q$ and $p\ga_p>2+N$ when $p>4+\f{4}{N}$, $Q_\mu(s\s u)=0$ has only one solution $s_\mu(u)\in\R$.

\item[(2)] From $(1)$, $Q_\mu(s\s u)>0$ when $s<s_\mu(u)$ and $Q_\mu(s\s u)<0$ when $s>s_\mu(u)$. So $I_\mu(s\s u)$ is 
	strictly increasing in $s\in(-\iy,s_\mu(u))$ and is strictly decreasing in $s\in(s_\mu(u),+\iy)$. Obviously, 
		$$\lim_{s\to-\iy}I_\mu(s\s u)=0^+,\quad \lim_{s\to+\iy}I_\mu(s\s u)=-\iy,$$
	which implies that 
		$$I_\mu(s_\mu(u)\s u)=\max_{s\in\R}I_\mu(s\s u)>0.$$

\item[(3)] It can be obtained directly from $(2)$.

\item[(4)] Let $\Phi_\mu(s,u)=Q_\mu(s\s u)$. Then $\Phi_\mu(s_\mu(u), u)=0$. Moreover,
\be
	\begin{aligned}
	\f{\pa}{\pa s} \Phi_\mu(s,u)&=\q(1+\ga_\q)^2\mu e^{\q(1+\ga_\q)s}\int_\RN\gu{u}\q+2e^{2s}\int_\RN\gu{u}2\\
			&\quad~~ +(2+N)^2e^{(2+N)s}\int_\RN|u|^2\gu{u}2-p\ga_p^2e^{p\ga_ps}\int_\RN|u|^p.
	\end{aligned}
\ee
Combining with $Q_\mu(s_\mu(u)\s u)=0$, we obtain
\be
	\begin{aligned}
	\f{\pa}{\pa s} \Phi_\mu(s_\mu(u),u)&=-(p\ga_p-\q-\q\ga_\q)(1+\ga_\q)\mu\int_\RN\gu{u}\q-(p\ga_p-2)\int_\RN\gu{u}2\\
		&\quad~~ -(p\ga_p-2-N)(2+N)\int_\RN|u|^2\gu{u}2<0.
	\end{aligned}
\ee
Then the Implict Function Theorem \cite{Chang=2005} implies that the map $u\mapsto s_\mu(u)$ is of class $\CR^1$.

\item[(5)] Since
	$$Q_\mu(s_\mu(u)\s (-u))=Q_\mu(-s_\mu(u)\s u)=Q_\mu(s_\mu(u)\s u)=0,$$
by the uniqueness, there is $s_\mu(-u)=s_\mu(u)$.
\end{itemize}
\ep

\subsection{Ground state critical point of \texorpdfstring{$I_\mu|_{\SR(a)}$}{} }\lab{groundstatesolution}

	In this subsection, we consider a minimizition problem
\be\lab{min}
	m_\mu(a):=\inf_{u\in\QR_\mu(a)} I_\mu(u).
\ee
From Lemma \ref{lempho}, we know that if $m_\mu(a)$ is achieved, then the minimizer is a ground 
state critical point of $I_\mu|_{\SR(a)}$.
We have the following lemma.
\bl\lab{est}
	The following statements hold.
	\begin{itemize}[fullwidth,itemindent=0em]
	\item[(1)]	$\DR(a):=\inf_{0<\mu\le1,u\in\QR_\mu(a)}\int_\RN|u|^2\gu{u}2>0$ is independent of $\mu$.
	\item[(2)]	If $\sup_{n\ge1}I_\mu(u_n)<+\iy$ for $u_n\in\QR_\mu(a)$, then
			$$\sup_{n\ge1} \max \lbr{ \mu\int_\RN\gu{u}\q, \int_\RN|u|^2\gu{u}2, \int_\RN\gu{u}2}<+\iy.$$
	\end{itemize}
\el
\bp
\begin{itemize}[fullwidth,itemindent=0em]
\item[(1)]  For any $u\in\QR_\mu(a)$, by the inequality \eqref{GNine}, there holds
\be
	(2+N)\int_\RN|u|^2\gu{u}2\le\ga_p\int_\RN|u|^p\le 
			K(p,N)\ga_p a^{\f{4N-p(N-2)}{2(N+2)}}\sbr{\int_\RN |u|^2\gu{u}2}^{\f{N(p-2)}{2(N+2)}}.
\ee
Since $\f{N(p-2)}{2(N+2)}>1$, we obtain $\DR(a)>0$.

\item[(2)]  For any $u\in\QR_\mu(a)$, there is
\be\lab{minus}
	\begin{aligned}
	I_\mu(u)&=I_\mu(u)-\f{1}{p\ga_p}Q_\mu(u)\\
		&=\f{p\ga_p-\q-\q\ga_\q}{\q p\ga_p}\mu\int_\RN\gu{u}\q+\f{p\ga_p-2}{2p\ga_p}\int_\RN\gu{u}2
			+\f{p\ga_p-2-N}{p\ga_p}\int_\RN|u|^2\gu{u}2.
	\end{aligned}
\ee
So the conclusion holds.
\end{itemize}
\ep
\br\lab{estofmin}
	Form \eqref{minus}, we see that
		$$m_\mu(a)\ge \DR_0(a):= \f{p\ga_p-2-N}{p\ga_p} \DR(a)>0,\quad  \forall\mu\in(0,1].$$
\er
Then we have the following result.
\bl\lab{tem4}
	There exists a small $\rho>0$ independent of $\mu$ such that for any $0<\mu\le1$, we have that 
		$$0<\sup_{u\in B_\mu(\rho,a)}I_\mu(u)<\DR_0(a)\quad \text{and}\quad I_\mu(u),Q_\mu(u)>0,~~\forall u\in B_\mu(\rho,a),$$
where 
	$$B_\mu(\rho,a)=\lbr{u\in\SR(a):\mu\int_\RN\gu{u}\q+\int_\RN\gu{u}2+\int_\RN|u|^2\gu{u}2\le\rho}.$$
\el
\bp
From the definition of $I_\mu$, we have
	$$\sup_{u\in B_\mu(\rho,a)} I_\mu(u)\le\max\lbr{\f{1}{\q},\f{1}{2},1}\rho<\DR_0(a),$$
where $\rho>0$ is small and is independent of $\mu$. On the other hand, by inequality \eqref{GNine}, for any $u\in \pa B_\mu(r,a)$ 
with $0<r<\rho$ for a smaller $\rho>0$,
\be
	\begin{aligned}
	\inf_{\pa B_\mu(r,a)}I_\mu(u)&\ge \f{\mu}{\q}\int_\RN\gu{u}\q+\f{1}{2}\int_\RN\gu{u}2+\int_\RN|u|^2\gu{u}2\\
				&\quad~~-\f{K(p,N)}{p}a^{\f{4N-p(N-2)}{2(N+2)}}\sbr{\int_\RN |u|^2\gu{u}2}^{\f{N(p-2)}{2(N+2)}}\\
				&\ge\f{\mu}{\q}\int_\RN\gu{u}\q+\f{1}{2}\int_\RN\gu{u}2+C\int_\RN|u|^2\gu{u}2\\
				&\ge C_1(a,\q,p,N)r>0,\\
	\inf_{\pa B_\mu(r,a)}Q_\mu(u)&\ge C_2(a,\q,p,N)r>0,
	\end{aligned}
\ee
which finish the proof.
\ep

	To find a Palais-Smale sequence, we consider an auxilary funtional as the one in \cite{Jeanjean=1997},
\be\lab{auxfun}
	J_\mu(s,u):=I_\mu(s\s u):\R\times\XR\to\R.
\ee
We study $J_\mu$ on the radial space $\R\times\SR_r(a)$ with
	$$\SR_r(a):=\SR(a)\cap\XR_r,\quad \XR_r=W^{1,\q}_{rad}(\RN)\cap W^{1,2}_{rad}(\RN).$$
Notice that $J_\mu$ is of class $\CR^1$. By the Symmetric Critical Point Principle \cite{Palais=1979}, 
a Palais-Smale sequence for $J_\mu|_{\R\times\SR_r(a)}$ is also a Palais-Smale sequence for $J_\mu|_{\R\times\SR(a)}$.
Denoting the closed sublevel set by 
	$$I_\mu^c=\lbr{u\in\SR(a):I_\mu(u)\le c},$$ 
we introduce the minimax class
\be\lab{path}
	\Ga_\mu:=\lbr{\ga=(\al,\beta)\in\CR([0,1],\R\times\SR_r(a)):\ga(0)\in\{0\}\times B_\mu(\rho,a),\ga(1)\in\{0\}\times I_\mu^0},
\ee
with the associated minimax level
\be\lab{level}
	\sigma_\mu(a):=\inf_{\ga_\in\Ga_\mu}\sup_{t\in[0,1]}J_\mu(\ga(t)).
\ee
Then
\bl\lab{mountain-pass}
	For any $0<\mu\le1$, $m_\mu(a)=\sigma_\mu(a)$.
\el
\bp
	For any $\ga=(\al,\beta)\in\Ga_\mu$, let us consider the function
	$$f_\ga(t):=Q_\mu(\al(t)\s\beta(t)).$$
We have $f_\ga(0)=Q_\mu(\beta(0))>0$, by Lemma \ref{tem4}. We {\it claim}   that $f_\ga(1)=Q_\mu(\beta(1))<0$: indeed, since $I_\mu(\beta(1))<0$, 
we have that $s_\mu(\beta(1))<0$, which means that $Q_\mu(\beta(1))<0$ by Lemma \ref{structure}. Moreover, $f_\ga$ is continuous, and hence 
we deduce that there exists $t_\ga\in(0,1)$ such that $f_\ga(t_\ga)=0$, namely $\al(t_\ga)\s\beta(t_\ga)\in\QR_\mu(a)$. So
	$$\max_{t\in[0,1]}J_\mu(\ga(t))\ge I_\mu(\al(t_\ga)\s\beta(t_\ga))\ge m_\mu(a),$$
and consequently $\sigma_\mu(a)\ge m_\mu(a)$. 

	On the other hand, if $u\in\QR_\mu(a)\cap\XR_r$, then
	$$\ga_u(t):=(0,((1-t)s_0+ts_1)\s u)\in \Ga_\mu,$$
where $s_0\ll-1$ and $s_1\gg1$. Since
	$$I_\mu(u)\ge\max_{t\in[0,1]}I_\mu(((1-t)s_0+ts_1)\s u)\ge\sigma_\mu(a),$$
there holds 
	$$m_\mu^r(a):=\inf_{u\in\QR_\mu(a)\cap\XR_r}I_\mu(u)\ge\sigma_\mu(a).$$
Finally the inequality $m_\mu(a)\ge m_\mu^r(a)$ can be obtained easily by using the Symmetric decreasing rearrangement, see \cite{Loss-Lieb}.
\ep
\br\lab{monotonicity}
	For any $0<\mu_1<\mu_2\le1$, since $I_{\mu_2}(u)\ge I_{\mu_1}(u)$ and $\Ga_{\mu_2}\subset\Ga_{\mu_1}$, there holds
	$$\sigma_{\mu_2}(a)=\inf_{\ga_\in\Ga_{\mu_2}}\sup_{t\in[0,1]}J_{\mu_2}(\ga(t))
		\ge\inf_{\ga_\in\Ga_{\mu_2}}\sup_{t\in[0,1]}J_{\mu_1}(\ga(t))
		\ge\inf_{\ga_\in\Ga_{\mu_1}}\sup_{t\in[0,1]}J_{\mu_1}(\ga(t))=\sigma_{\mu_1}(a),$$
i.e.,  $\sigma_\mu(a)$ is non-decreasing with respect to $\mu\in(0,1]$.
\er

	We recall the following definition and theorem from \cite{Ghoussoub=1993}.

\newtheorem*{defA}{Definition A. \cite[Definition 3.1]{Ghoussoub=1993}}
\newtheorem*{thmB}{Theorem B. \cite[Theorem 5.2]{Ghoussoub=1993}}

\begin{defA}\lab{defA}
	Let $B$ be a closed subset of $X$. We say that a class $\FR$ of compact subsets of $X$ is a homotopy stable family with boundary
	 $B$ provided
	\begin{itemize}[fullwidth,itemindent=0em]
	\item[(a)]  every set in $\FR$ contains $B$.
	\item[(b)]  for any set $A$ in $\FR$ and any $\eta\in\CR([0,1]\times X,X)$ satisfying $\eta(t,x)=x~$ for all $(t,x)$ in 
		$(\{0\}\times X)\cup([0,1]\times B)$ we have that $\eta(1,A)\subset\FR$.
	\end{itemize}
\end{defA}
	
	We remark that the case $B=\emptyset$ is admissible.

\begin{thmB}\lab{thmB}
	Let $\phi$ be a $\CR^1$-functional on a complete connected $C^1$-Finsler manifold $X$ and  consider a homotopy stable family $\FR$
	with an extended closed boundary $B$. Set $c=c(\phi,\FR)$ and let $F$ be a closed subset of $X$ satisfying
	\be\lab{tem5}
		A\cap F\setminus B\neq\emptyset\quad \forall A\in\FR
	\ee 
	and
	\be\lab{tem6}
		\sup\phi(B)\le c\le \inf\phi(F).
	\ee
	Then for any sequence of sets $A_n\subset\FR$ such that $\lim_{n\to\iy}\sup_{A_n}\phi=c$, there exists a sequence
	 $x_n\subset X\setminus B$  such that
	\begin{itemize}[fullwidth,itemindent=0em]
	\item[(1)]  $\lim_{n\to\iy}\phi(x_n)=c$,
	\item[(2)]  $\lim_{n\to\iy}\nm{\rd\phi(x_n)}=0$,
	\item[(3)]  $\lim_{n\to\iy}\text{dist}(x_n,F)=0$,
	\item[(4)]  $\lim_{n\to\iy}\text{dist}(x_n,A_n)=0$.
	\end{itemize}
\end{thmB}

	Now we establish a technical result showing the existence of a Palais-Smale sequence of $\sigma_\mu(a)$ with an additional 
property.
\bl\lab{PSseq}
	For any fixed $\mu\in(0,1]$, there exists a sequence $u_n\in\SR_r(a)$ such that
		$$I_\mu(u_n)\to\sigma_\mu(a),\quad I_\mu|_{\SR(a)}'(u_n)\to0,\quad Q_\mu(u_n)\to0\quad\text{and}\quad u_n^-\to0 \text{ a.e. in }\RN.$$
\el
\bp
	Using Definition \ref{defA}, it is easy to check that $\FR=\lbr{A=\ga([0,1]):\ga\in\Ga_\mu}$ is a homotopy stable family of compact 
subsets of $X=\R\times\SR_\mu^r$ with boundary $B=(\{0\}\times B_\mu(\rho,a))\cup(\{0\}\times I_\mu^0)$. Set
 $F=\lbr{J_\mu\ge\sigma_\mu(a)}$, then the assumptions \eqref{tem5} and \eqref{tem6} with $\phi=J_\mu$, $c=\sigma_\mu(a)$ are 
satisfied. Therefore, taking a minimizing 
sequence $\lbr{\ga_n=(0,\beta_n)}\subset\Ga_\mu$ with $\beta_n\ge0$ a.e. in $\RN$, there exists a Palais-Smale sequence
 $\lbr{(s_n,w_n)}\subset\R\times\SR_r(a)$ for $J_\mu|_{\R\times\SR_r(a)}$ at level $\sigma_\mu(a)$, that is
\be\lab{tem8}
	\pa_s J_\mu(s_n,w_n)\to0\quad\text{and}\quad \pa_u J_\mu(s_n,w_n)\to0\quad \text{as }n\to\iy,
\ee
with the additional property that
\be\lab{tem7}
	|s_n|+\text{dist}_\XR(w_n,\beta_n([0,1]))\to0 \quad \text{as }n\to\iy.
\ee
Let $u_n=s_n\s w_n$. 
The first condition in \eqref{tem8} reads $Q_\mu(u_n)\to0$, while the second condition gives
\be
	\begin{aligned}
	\nm{\rd I_\mu|_{\SR(a)}(u_n)}&=\sup_{\psi\in T_{u_n}\SR(a),\nm{\psi}_\XR\le1}|\rd I_\mu(u_n)[\psi]|\\
		&=\sup_{\psi\in T_{u_n}\SR(a),\nm{\psi}_\XR\le1}|\rd I_\mu(s_n\s w_n)[s_n\s (-s_n)\s\psi]|\\
		&=\sup_{\psi\in T_{u_n}\SR(a),\nm{\psi}_\XR\le1}|\pa_u J_\mu(s_n,w_n)[(-s_n)\s\psi]|\\
		&\le \nm{\pa_u J_\mu(s_n,w_n)} \sup_{\psi\in T_{u_n}\SR(a),\nm{\psi}_\XR\le1}|(-s_n)\s\psi|\\
		&\le C \nm{\pa_u J_\mu(s_n,w_n)}\to0\quad \text{as }n\to\iy.
	\end{aligned}
\ee
Finally, \eqref{tem7} implies that $u_n^-\to0$ a.e. in $\RN$.
\ep

	Now we show the compactness of the Palais-Smale sequence obtained in Lemma \ref{PSseq}.
\bl\lab{PScon}
	For any fixed $\mu\in(0,1]$, let $u_n$ be a sequence obtained in Lemma \ref{PSseq}. Then there exists  a
 $u_\mu\in\XR\setminus\{0\}$ and a $\la_\mu\in\R$ such that 
up to a subsequence,
\be
	u_n\rh u_\mu\ge0\quad \text{in }\XR,
\ee
\be
	I_\mu(u_\mu)=\sigma_\mu(a)\quad\text{and}\quad I_\mu'(u_\mu)+\la_\mu u_\mu=0.
\ee
Moreover, if $\la_\mu\neq0$, we have that
	$$u_n\ra u_\mu\quad \text{in }\XR.$$
\el
\bp
	From Lemma \ref{est} and Remark \ref{monotonicity}, we know that $u_n$ is bounded in $\XR_r$. Thus by 
\cite[Propositon 1.7.1]{Cazenave}, we conclude that up to a subsequence, there exists a $u_\mu\in\XR_r$ such that
	$$u_n\rh u_\mu\quad \text{in }\XR \text{ and in } L^2(\RN),$$
	$$u_n\ra u_\mu\quad \text{in }L^q(\RN),~\forall q\in(2,2^*),$$
	$$u_n\ra u_\mu\ge0\quad \text{a.e. in }\R.$$
By interpolation and inequality \eqref{GNine}, we have that
	$$u_n\ra u_\mu\quad \text{in }L^q(\RN),~\forall q\in(2,22^*).$$
We claim that $u_\mu\neq0$. Assume $u_\mu=0$, then as $n\to\iy$
	$$(1+\ga_\q)\mu\int_\RN|\nabla u_n|^\q+\int_\RN|\nabla u_n|^2+(2+N)\int_\RN|u_n|^2|\nabla u_n|^2=Q_\mu(u_n)+\ga_p\int_\RN|u_n|^p\to0, $$
which implies that $I_\mu(u_n)\to0$, in contradiction with Remark \ref{estofmin}. So $u_\mu\neq0$.
By \cite[Lemma 3]{Berestycki-Lions=1983-2}, it follows from  $I_\mu|_{\SR(a)}'(u_n)\to0$ that there exists a 
sequence $\la_n\in\R$ such that
\be\lab{tem10}
	I_\mu'(u_n)+\la_n u_n\to0\quad \text{in}~\XR^*.
\ee
Hence $\la_n=\f{1}{a}I_\mu'(u_n)[u_n]+o_n(1)$ is bounded in $\R$, and we assume, up to a subsequence, $\la_n\to\la_\mu$. Since $u_n$ is 
bounded, we have $I_\mu'(u_n)+\la_\mu u_n\to0$. From Lemma \ref{A2}, we see that 
\be\lab{tem9}
	I_\mu'(u_\mu)+\la_\mu u_\mu=0.
\ee
Then testing \eqref{tem9} with $x\cdot\nabla u$ and $u$, we obtain $Q_\mu(u_\mu)=0$. It follows that 
	$$Q_\mu(u_n)+\ga_p\int_\RN|u_n|^p\to Q_\mu(u_\mu)+\ga_p\int_\RN|u_\mu|^p.$$
Then using the weak lower semicontinuous property, see \cite[Lemma 4.3]{Colin-Jeanjean-Squassina=2010}, there must be
\be\lab{tem11}
	\mu\int_\RN|\nabla u_n|^\q\to \mu\int_\RN|\nabla u_\mu|^\q,
\ee
\be\lab{tem12}
	\int_\RN|\nabla u_n|^2\to\int_\RN|\nabla u_\mu|^2,
\ee
\be\lab{tem13}
	\int_\RN|u_n|^2|\nabla u_n|^2\to\int_\RN|u_\mu|^2|\nabla u_\mu|^2.
\ee
That gives $I_\mu(u_\mu)=\lim_{n\to\iy}I_\mu(u_n)=\sigma_\mu(a)$. Moreover, from \eqref{tem11}-\eqref{tem13}, we obtain
\be\lab{tem14}
	I_\mu'(u_n)[u_n]\to I_\mu'(u_\mu)[u_\mu].
\ee
Thus combining \eqref{tem14} with \eqref{tem10}-\eqref{tem9}, there holds $\la_\mu\nm{u_n}_2^2\to\la_\mu\nm{u_\mu}_2^2$. So $\la_\mu\neq0$ 
implies that $u_n\ra u_\mu$ in $\XR$.
\ep

	Based on the above preliminary works, we conclude that
\bt\lab{groundstate}
	For any fixed $\mu\in(0,1]$, there exists a $u_\mu\in\XR_r\setminus\{0\}$ and a $\la_\mu\in\R$ such that 
\begin{align*}
	&I_\mu'(u_\mu)+\la_\mu u_\mu=0,\\
	&I_\mu(u_\mu)=m_\mu(a),\quad Q_\mu(u_\mu)=0,\\
	&0<\nm{u_\mu}_2^2\le a,\quad u_\mu\ge0.
\end{align*}
Moreover, if $\la_\mu\neq0$, we have that $\nm{u_\mu}_2^2=a$, i.e.,  $m_\mu(a)$ is achieved, and $u_\mu$ is 
a ground state critical point of $I_\mu|_{\SR(a)}$.
\et

\bp[Proof of Theorem \ref{thm1} for the case $N=1$:]
	When $N=1$, there is $W^{1,2}(\R)\hra\CR^{0,\al}(\R)$, so $V(u)$ and hence $I(u)$ is of class $\CR^1(W^{1,2}(\R))$.
Then one can follow  the process in this subsection to prove Theorem \ref{thm1} by taking $\mu=0$, 
but we claim that there needs some necessary modifications,
since the compact embedding  $W^{1,2}_{rad}(\RN)\hra\hra L^q(\RN)$ for $2<q<2^*$ does not hold when $N=1$. However 
the compactness still holds for bounded sequences of radially decreasing functions (see e.g. \cite[Propositon 1.7.1]{Cazenave}).
So we need to confirm that the Palais-Smale sequence obtained in Lemma \ref{PSseq} consists of radially decreasing functions. 
Then it is natural to replace the minimizing sequence $\ga_n=(0,\beta_n)$ choosen in Lemma \ref{PSseq} 
with $\bar \ga_n:=(0,\bar \beta_n)$, where $\bar \beta_n(t)=|\beta_n(t)|^*$ is the symmetric decreasing rearrangement of $\beta_n(t)$ 
at every $t\in[0,1]$.
This is a natural candidate to be minimizing sequence, with $\bar \beta_n(t)\ge0$, radially symmetric and decreasing 
for every $t\in[0,1]$. In order to check that $\bar \ga_n\in\Ga_0$, we have to check that each $\bar\beta_n$ is continuous 
on $[0,1]$, which has been proved in \cite{Coron=1984} (for more argument we refer to \cite[Remark 5.2]{Soave=JDE=2020}). 
As a result, Theorem \ref{groundstate} with $\mu=0$ holds, and combining with Lemma \ref{Liouville1}, we obtain 
Theorem \ref{thm1} immediately.
\ep

\subsection{Infinitely many critical points of \texorpdfstring{$I_\mu|_{\SR(a)}$}{} }\lab{infinitelymanysolutions}

	This subsection concerns the existence of infinitely many radial critical points of $I_\mu|_{\SR(a)}$.
 Denote $\tau(u)=-u$ and let $Y\subset\XR$. 
A set $A\subset Y$ is called $\tau$-invariant if $\tau(A)=A$. A homotopy $\eta:[0,1]\times Y\to Y$ is  $\tau$-equivariant if 
$\eta(t,\tau(u))=\tau(\eta(t,u))$ for all $(t,u)\in[0,1]\times Y$. We recall the following definition.

\newtheorem*{defC}{Definition C. \cite[Definition 7.1]{Ghoussoub=1993}}

\begin{defC}\lab{defC}
	Let $B$ be a closed $\tau$-invariant subset of $Y$. A class $\GR$ of compact subsets of $Y$ is said to be
	a $\tau$-homotopy stable family with boundary $B$ provided
	\begin{itemize}[fullwidth,itemindent=0em]
	\item[(a)]	every set in $\GR$ is $\tau$-invariant,
	\item[(b)]  every set in $\GR$ contains $B$,
	\item[(c)]  for any set $A\in\GR$ and any $\tau$-equivariant homotopy $\eta\in\CR([0,1]\times Y,Y)$ satisfying $\eta(t,x)=x$ 
		for all $(t,x)$ in $(\{0\}\times Y)\cup([0,1]\times B)$ we have that $\eta(1,A)\subset\GR$.
	\end{itemize}
\end{defC}

	Following the strategy of \cite[Section. 5]{Jeanjean-Lu=CVPDE=2020}, we consider 
the functional $K_\mu:\XR\setminus\{0\}\to\R$ defined by
\be\lab{auxilary1}
	\begin{aligned}
	K_\mu(u):&=I_\mu(s_\mu(u)\s u)\\
			&=\f{\mu}{\q} e^{\q(1+\ga_\q)s_\mu(u)}\int_\RN\gu{u}\q+\f{1}{2}e^{2s_\mu(u)}\int_\RN\gu{u}2\\
			&\quad~~ +e^{(2+N)s_\mu(u)}\int_\RN|u|^2\gu{u}2-\f{1}{p} e^{p\ga_ps_\mu(u)}\int_\RN|u|^p,
	\end{aligned}
\ee
where $s_\mu(u)$ is given by Lemma \ref{structure}. Then we  see that $K_\mu(u)$ is $\tau$-invariant. Moreover, inspired by 
\cite[Proposition 2.9]{Szulkin-Weth=JFA=2009}, there holds

\bl\lab{equiv}
	The functional $K_\mu$ is of class $\CR^1$ and 
	\begin{align*}
		K_\mu'(u)[\phi]&=\mu e^{\q(1+\ga_\q)s_\mu(u)}\int_\RN\gu{u}{\q-2}\nabla u\cdot\nabla\phi  
							+e^{2s_\mu(u)}\int_\RN\nabla u\cdot\nabla\phi\\
				&\quad~~ +2e^{(2+N)s_\mu(u)}\int_\RN(u\phi\gu{u}2+|u|^2\nabla u\cdot\nabla\phi) 
							- e^{p\ga_ps_\mu(u)}\int_\RN|u|^{p-2}u\phi\\
				&=I_\mu'(s_\mu(u)\s u)[s_\mu(u)\s \phi],
	\end{align*}
	for any $u\in\XR\setminus\{0\}$ and $\phi\in\XR$.
\el
\bp
	Let $u\in\XR\setminus\{0\}$ and $\phi\in\XR$. We estimate the term
		$$K_\mu(u_t)-K_\mu(u)=I_\mu(s_t\s u_t)-I_\mu(s_0\s u),$$
where $u_t=u+t\phi$ and $s_t=s_\mu(u_t)$ with $|t|$ small enough. By the mean value theorem, we have
\begin{align*}
	&\quad~~I_\mu(s_t\s u_t)-I_\mu(s_0\s u)\le I_\mu(s_t\s u_t)-I_\mu(s_t\s u)\\
	&=\mu e^{\q(1+\ga_\q)s_t}\int_\RN|\nabla u_{\eta_t}|^{\q-2}(\nabla u\cdot\nabla\phi+\eta_t|\nabla\phi|^2)t  
							+e^{2s_t}\int_\RN(\nabla u\cdot\nabla\phi+\f{t}{2}|\nabla\phi|^2)t\\
	&\quad~ +2e^{(2+N)s_t}\int_\RN\sbr{ u_{\eta_t}\phi|\nabla u_{\eta_t}|^2  
		       +|u_{\eta_t}|^2(\nabla u\cdot\nabla\phi+\eta_t|\nabla\phi|^2)}t  
							- e^{p\ga_ps_t}\int_\RN|u_{\eta_t}|^{p-2}(u\phi+\f{\eta_t}{2}\phi^2)t,
\end{align*}
where $|\eta_t|\in(0,|t|)$. Similarly, 
\begin{align*}
	&\quad~~I_\mu(s_t\s u_t)-I_\mu(s_0\s u)\ge I_\mu(s_0\s u_t)-I_\mu(s_0\s u)\\
	&=\mu e^{\q(1+\ga_\q)s_0}\int_\RN|\nabla u_{\xi_t}|^{\q-2}(\nabla u\cdot\nabla\phi+\xi_t|\nabla\phi|^2)t  
							+e^{2s_0}\int_\RN(\nabla u\cdot\nabla\phi+\f{t}{2}|\nabla\phi|^2)t\\
	&\quad~ +2e^{(2+N)s_0}\int_\RN\sbr{ u_{\xi_t}\phi|\nabla u_{\xi_t}|^2  
		       +|u_{\xi_t}|^2(\nabla u\cdot\nabla\phi+\xi_t|\nabla\phi|^2)}t  
							- e^{p\ga_ps_0}\int_\RN|u_{\xi_t}|^{p-2}(u\phi+\f{\xi_t}{2}\phi^2)t,
\end{align*}
where $|\xi_t|\in(0,|t|)$. Since $s_t\to s_0$ as $t\to0$, it follows from the two inequalities above that
\begin{align*}
	\lim_{t\to0}\f{K_\mu(u_t)-K_\mu(u)}{t}&=\mu e^{\q(1+\ga_\q)s_\mu(u)}\int_\RN\gu{u}{\q-2}\nabla u\cdot\nabla\phi  
							+e^{2s_\mu(u)}\int_\RN\nabla u\cdot\nabla\phi\\
			&\quad~~ +2e^{(2+N)s_\mu(u)}\int_\RN(u\phi\gu{u}2+|u|^2\nabla u\cdot\nabla\phi) 
							- e^{p\ga_ps_\mu(u)}\int_\RN|u|^{p-2}u\phi.
\end{align*}
Then similarly as Lemma \ref{A1}, we see that the G\^ateaux derivative of $K_\mu$ is bounded linear and continuous. 
Therefore $K_\mu$ is of class $\CR^1$, see \cite{Chang=2005}. In particular, by changing variables in the integrals, we have
	$$K_\mu'(u)[\phi]=I_\mu'(s_\mu(u)\s u)[s_\mu(u)\s \phi].$$
The proof is complete.
\ep

	To get  the particular Palais-Smale sequence of $I_\mu|_{\SR(a)}$ as the one in Lemma \ref{PSseq}, we need
\bl\lab{minmax}
	Let $\GR$ be a $\tau$-homotopy stable family of compact subsets of $Y=\SR_r(a)$ with boundary $B=\emptyset$, 
	and set
		$$d:=\inf_{A\in\GR}\max_{u\in A}K_\mu(u).$$
	If $d>0$, then there exists a sequence $u_n\in \SR_r(a)$ such that
	$$I_\mu(u_n)\to d,\quad I_\mu|_{\SR(a)}'(u_n)\to0, \quad Q_\mu(u_n)=0.$$
\el
\bp
	Let $A_n\in\GR$ be a minimizing sequence of $d$. We define the mapping 
	$$\eta:[0,1]\times \SR(a)\to\SR(a),\quad \eta(t,u)=(ts_\mu(u))\s u,$$
which is continuous and satisfies $\eta(t,u)=u$ for all $(t,u)\in\{0\}\times\SR(a)$. Thus, by the definition of $\GR$, 
one has
	$$D_n:=\eta(1,A_n)=\lbr{s_\mu(u)\s u:u\in A_n}\in\GR.$$
In particular, $D_n\subset\QR_\mu(a)$ for any $n\in\N^+$. For any $u\in\SR(a)$ and any $s\in\R$, we see that
	$$Q_\mu((s_\mu(u)-s)\s (s\s u))=Q_\mu(s_\mu(u)\s u))=0,$$
that is $s_\mu(s\s u)=s_\mu(u)-s$, which gives $K_\mu(s\s u)=K_\mu(u)$. Then it is clear that 
$\max_{D_n}K_\mu=\max_{A_n}K_\mu\to d$ and thus $D_n$ is another minimizing sequence of $d$. 
Now, using the minimax principle \cite[Theorem 7.2]{Ghoussoub=1993}, we  obtain a Palais-Smale sequence $v_n\in\SR(a)$ 
for $K_\mu$ at the level $d$ such that 
	$$\text{dist}_\XR(v_n, D_n)\to0.$$
Finally, a similar argument as the one in Lemma \ref{PSseq} gives that $u_n=s_n\s v_n$ satisfying that
	$$I_\mu(u_n)\to d,\quad I_\mu|_{\SR(a)}'(u_n)\to0, \quad Q_\mu(u_n)=0.$$
\ep

	To construct a sequence of $\tau$-homotopy stable families of compact subsets of $\SR_r(a)$ with boundary $B=\emptyset$, 
we proceed as in \cite[Section. 8]{Berestycki-Lions=1983-2}. Since $\XR$ is separable, there exists a nested sequence of finite 
dimensional subspaces of $\XR$, $W_1\subset W_2\subset \cdots\subset W_i\subset W_{i+1}\subset\cdots\subset \XR$ such that 
$dim(W_i)=i$ and the closure of $\cup_{i\in\N^+} W_i$ in $\XR$ is equal to $\XR$. Note that since $\XR$ is dense in $W^{1,2}(\RN)$,
the closure in $W^{1,2}(\RN)$ is also equal to $W^{1,2}(\RN)$. Since $W^{1,2}(\RN)$ is a Hilbert space, 
we denote by $P_i$ the orthogonal projection from $W^{1,2}(\RN)$ onto $W_i$. We also recall the definition of the genus of 
$\tau$-invariant sets due to M. A. Krasnoselskii and refer to \cite[Section. 7]{Rabinowitz=1986}.

\newtheorem*{defD}{Definition D. (Krasnoselskii genus)}

\begin{defD}\lab{defD}
	For any nonempty closed $\tau$-invariant set $A\subset\XR$. The genus of $A$ is defined by
		$$\ind(A):=\min\lbr{ k\in\N^+:\exists~\phi:A\to\R^k\setminus\{0\}, \phi \text{ is odd and continuous }}.$$
	We set $\ind(A)=+\iy$ if such $\phi$ does not exist, and set $\ind(A)=0$ if $A=\emptyset$.
\end{defD}

	Let $\AR(a)$ be the family of compact $\tau$-invariant subsets of $\SR_r(a)$. For each $j\in\N^+$, set 
	$$\AR_j(a):=\lbr{A\in\AR(a):\ind(A)\ge j}$$
and 
	$$c_\mu^j(a):=\inf_{A\in\AR_j(a)}\max_{u\in A}K_\mu(u).$$
Concerning $\AR_j(a)$ and $c_\mu^j(a)$, we have
\bl\lab{genuslevel}
	\begin{itemize}[fullwidth,itemindent=0em]
	\item[(1)]	For any $j\in\N^+$, $\AR_j(a)\neq\emptyset$ and $\AR_j(a)$ is a $\tau$-homotopy stable family of 
		compact subsets of $\SR_r(a)$ with boundary $B=\emptyset$.
	\item[(2)]	For any $\mu\in(0,1]$, any $j\in\N^+$, $c_\mu^{j+1}(a)\ge c_\mu^j(a)\ge\DR_0(a)>0$.
	\item[(3)]	For any $j\in\N^+$, $c_\mu^j(a)$ is non-decreasing with respect to $\mu\in(0,1]$.
	\item[(4)]	$b_j(a):=\inf_{0<\mu\le1}c_\mu^j(a)\to+\iy$ as $j\to+\iy$.
	\end{itemize}
\el
\bp
\begin{itemize}[fullwidth,itemindent=0em]
\item[(1)]	For any $j\in\N^+$, $\SR_r(a)\cap W_j\in\AR(a)$. By the basic properties of the genus, one has
	$$\ind(\SR_r(a)\cap W_j)=j$$
and thus $\AR_j(a)\neq\emptyset$. The rest is clear by the properties of the genus. 

\item[(2)]	For any $A \in\AR_j(a)$, using the fact that $s_\mu(u)\s u\in\QR_\mu(a)$ for all $u\in A$, we have
	$$\max_{u\in A}K_\mu(u)=\max_{u\in A}I_\mu(s_\mu(u)\s u)\ge m_\mu(a)\ge \DR_0(a)$$
and thus $c_\mu^j(a)\ge\DR_0(a)>0$. Since $\AR_{j+1}(a)\subset\AR_j(a)$, it is clear that $c_\mu^{j+1}(a)\ge c_\mu^j(a)$.

\item[(3)]	For any $0<\mu_1<\mu_2\le1$, any $u\in A\in\AR_j(a)$, there holds
	$$K_{\mu_2}(u)=I_{\mu_2}(s_{\mu_2}(u)\s u)\ge I_{\mu_2}(s_{\mu_1}(u)\s u)>I_{\mu_1}(s_{\mu_1}(u)\s u)=K_{\mu_1}(u),$$
which means $c_{\mu_2}^j(a)\ge c_{\mu_1}^j(a)$, i.e., $c_\mu^j(a)$ is non-decreasing with respect to $\mu\in(0,1]$.

\item[(4)]	The proof is inspired by that of \cite[Theorem 9]{Berestycki-Lions=1983-2}. 
First, we claim that

\vskip 0.05in
\noindent {\bf Claim: }{\it for any $M>0$, there exists a small $\delta_0=\delta_0(a,M)>0$, a small $r_0=r_0(a,M)>0$ and a large $k_0=k_0(a,M)\in\N^+$ 
	such that for any $0<\mu<\delta_0$ and any $k\ge k_0$, one has 
		$$I_\mu(u)\ge M\quad \text{if}\quad \nm{P_ku}_\XR\le r_0~\text{and}~u\in \QR_\mu^r(a).$$}
Now we check it. By contradiction, we assume that there exists $M_0>0$ such that for any $0<\delta\le1$, any $r>0$ and any $k\in\N^+$ 
one can always find $\mu\in(0,\delta]$, $l\ge k$ and $u\in\QR_\mu^r(a)$ such that
	$$\nm{P_ku}_\XR\le r \quad \text{but}\quad I_\mu(u)< M_0.$$
As a consequnce, one can obtain a sequence $\mu_n\to0^+$, a sequence $k_n\to+\iy$, and a sequence $u_n\in\QR_{\mu_n}^r(a)$ such that
	$$\Vert P_{k_n} u_n\Vert_\XR\le\f{1}{n}\quad \text{and}\quad I_{\mu_n}(u_n)< M_0$$
for any $n\in\N^+$. From Lemma \ref{est}, we know that $u_n$ is bounded in $W^{1,2}(\RN)$. Since $P_{k_n} u_n$ is also bounded in $\XR$, 
we assume that up to a subsequence
	$$u_n\rh u~\text{ in}~W^{1,2}(\RN)\quad\text{and}\quad P_{k_n} u_n\rh v~\text{ in}~\XR.$$
We show that $u=v$. Indeed, one also has $P_{k_n} u_n\rh v$ in $W^{1,2}(\RN)$ and 
\begin{align*}
	\nm{u-v}_{W^{1,2}(\RN)}^2&=\lim_{n\to\iy}\left<u_n-P_{k_n} u_n,u-v\right>_{W^{1,2}(\RN)}\\
		&=\lim_{n\to\iy}\left<u_n,u-v\right>_{W^{1,2}(\RN)}-\lim_{n\to\iy}\left<P_{k_n} u_n,u-v\right>_{W^{1,2}(\RN)}\\
		&=\left<u,u-v\right>_{W^{1,2}(\RN)}-\lim_{n\to\iy}\left<u_n,P_{k_n} u-P_{k_n} v\right>_{W^{1,2}(\RN)}\\
		&=\left<u,u-v\right>_{W^{1,2}(\RN)}-\left<u,u-v\right>_{W^{1,2}(\RN)}\\
		&=0,
\end{align*}
where we use the fact that $P_{k_n} u\to u$ and $P_{k_n} v\to v$ in $W^{1,2}(\RN)$. 
Therefore $u=v$ and $u\in\XR$. Since $\Vert P_{k_n} u_n\Vert_\XR\to0$, there must be $u=0$. Then combining the interpolation inequality
and the fact that $\sup_{n\in\N^+}\int_\RN|u_n|^2|\nabla u_n|^2<+\iy$, we obtain $\nm{u_n}_p\to0$. Futher $u_n\in\QR_{\mu_n}(a)$ gives that 
	$$\mu_n\int_\RN|\nabla u_n|^\q\to0,\quad \int_\RN|\nabla u_n|^2\to0,\quad \int_\RN|u_n|^2|\nabla u_n|^2\to0,$$
which is in contradiction with Lemma \ref{est}. So we prove the claim. 

	Then we can prove the conclusion $(4)$. By contradiction, we assume that
		$$\liminf_{j\to\iy} b_j<M\quad \text{for some } M>0.$$
Then there exist $\mu\in(0,\delta_0)$, $k>k_0$ such that $c_\mu^k(a)<M$. By the definition of $c_\mu^k(a)$, one can find $A\in\AR_k(a)$ 
 such that
 	$$\max_{u\in A}I_\mu(s_\mu(u)\s u)=\max_{u\in A}K_\mu(u)<M.$$
Since Lemma \ref{structure} implies that the mapping $\vp:A\to\QR_\mu^r(a)$ defined by $\vp(u)=s_\mu(u)\s u$ is odd and continuous, 
we have $\bar A:=\vp(A)\subset\QR_\mu^r(a)$, $\max_{u\in\bar A}I_\mu(u)<M$ and 
\be\lab{tem15}
	\ind(\bar A)\ge\ind(A)\ge k>k_0.
\ee
On the other hand, it follows from the claim that $\inf_{u\in\bar A}\Vert P_{k_0} u_n\Vert_\XR\ge r_0>0$. Setting
	$$\psi(u)=\f{P_{k_0}u}{\Vert P_{k_0} u_n\Vert_\XR}\quad \text{for any }u\in\bar A,$$
we obtain an odd continuous mapping $\psi:\bar A\to\psi(\bar A)\subset W_{k_0}\setminus\{0\}$ and thus 
	$$\ind(\bar A)\le\ind(\psi(\bar A))\le k_0,$$
which contradicts \eqref{tem15}. Therefore we have $b_j(a)\to+\iy$ as $j\to+\iy$.
\end{itemize}
\ep

	For any fixed $\mu\in(0,1]$ and any $j\in\N^+$, by Lemma \ref{minmax} and \ref{genuslevel}, 
one can find a sequence $u_n\in\SR_r(a)$ such that
	$$I_\mu(u_n)\to c_\mu^j(a),\quad I_\mu|_{\SR(a)}'(u_n)\to0,\quad Q_\mu(u_n)=0.$$
Then similar to Lemma \ref{PScon}, we have
\bl
 	There exists  a $u_\mu^j\in\XR\setminus\{0\}$ and a $\la_\mu^j\in\R$ such that up to a subsequence,
\be
	u_n^j\rh u_\mu^j\quad \text{in }\XR,
\ee
\be
	I_\mu(u_\mu^j)=c_\mu^j(a)\quad\text{and}\quad I_\mu'(u_\mu^j)+\la_\mu^j u_\mu^j=0.
\ee
Moreover, if $\la_\mu^j\neq0$, we have that
	$$u_n^j\ra u_\mu^j\quad \text{in }\XR.$$
\el

	Based on the above preliminary works, we conclude that
\bt\lab{infinite}
	For any fixed $\mu\in(0,1]$ and any $j\in\N^+$, there exists  a $u_\mu^j\in\XR_r\setminus\{0\}$ and a $\la_\mu^j\in\R$ such that 
\begin{align*}
	&I_\mu'(u_\mu^j)+\la_\mu^j u_\mu^j=0,\\
	&I_\mu(u_\mu^j)=c_\mu^j(a),\quad Q_\mu(u_\mu^j)=0,\\
	&0<\nm{u_\mu^j}_2^2\le a.
\end{align*}
Moreover, if $\la_\mu^j\neq0$, we have that $\nm{u_\mu^j}_2^2=a$, i.e.,  $\lbr{u_\mu^j:j\in\N^+}$ are infinitely many critical points
of $I_\mu|_{\SR(a)}$ with increasing energy.
\et

\vskip0.23in
\section{Convergence issues as \texorpdfstring{$\mu\to0^+$}{}}\lab{convergethm}
	
	In this section, letting $\mu\to0^+$, we show that the sequences of critical points of $I_\mu|_{\SR(a)}$ obtained in 
the Section \ref{sec3} converge to critical points of $I|_{\tilde\SR(a)}$.

\bt\lab{converge}
	Let $N\ge2$. Suppose that $\mu_n\to0^+$, $I_\n'(u_\n)+\la_\n u_\n=0$ with $\la_\n\ge0$ and $I_\n(u_\n)\to c\in(0,+\iy)$ 
for $u_\n\in\SR_r(a_n)$ with $0<a_n\le a$. 
Then there exists a subsequence $u_\n\rh u$ in $W^{1,2}(\RN)$ with $u\neq0$, $u\in W^{1,2}_{rad}(\RN)\cap L^\iy(\RN)$ and
there exists a $\la\in\R$ such that 
	$$\quad I'(u)+\la u=0,\quad I(u)=c \quad\text{and}\quad 0<\nm{u}_2^2\le a.$$
Moreover, 
\begin{itemize}[fullwidth,itemindent=0em]
\item[(1)]	if $u_\n\ge0$ for each $n\in\N^+$, then $u\ge0$,
\item[(2)]	if $\la\neq0$, we have that $\nm{u}_2^2=\lim_{n\to\iy}a_n$.
\end{itemize}
\et
\br\lab{remark2}
	We note that the condition $\la_\n\ge0$ is only used in the following Step 1 to realize the Morse iteration. 
	If one can prove the conclusion in Step 1 without this condition, then the conclusion in Theorem \ref{thm1} can 
	be extended to $N=3,4$ with $4+\f{4}{N}<p<22^*$.
\er
\bp[Proof of Theorem \ref{converge}:]
	The proof is inspired by \cite{Jeanjean-Luo-Wang=JDE=2015,MathN=2019}. 
First, by Lemma \ref{lempho}, $I_\n'(u_\n)+\la_\n u_\n=0$ implies that 
	$$Q_\n(u_\n)=0\quad \text{for each }n\in\N^+.$$ 
Then from Lemma \ref{est}, we see that
\be\lab{bdd}
	\sup_{n\ge1} \max \lbr{ \n\int_\RN|\nabla u_n|^\q, \int_\RN|u_n|^2|\nabla u_n|^2, \int_\RN|\nabla u_n|^2}<+\iy,
\ee
and hence $u_\n$ is bounded in $W^{1,2}(\RN)$. We claim that $\liminf_{n\to\iy}a_n>0$ and hence $\la_\n=\f{1}{a_n}I_\n'(u_\n)$
 $[u_\n]$ is also bounded in $\R$. Indeed, if $a_n\to0$, then $\nm{u_\n}_p\to0$, and it follows from
 $\QR_\n(u_n)=0$ that $I_\n(u_\n)\to0$ which contradicts that $c>0$.
Thus, up to a subsequence, $\la_\n\to\la$ in $\R$, $u_\n\rh u$ in $W^{1,2}_{rad}(\RN)$, $u_\n\to u$ in $L^q(\RN)$ for $2<q<22^*$, 
and $u_\n\to u$ a.e. on $\RN$. So if $u_\n\ge0$ for each $n\in\N^+$, we have that $u\ge0$. Moreover, a similar argument 
as Lemma \ref{A2} tells that $u_n\nabla u_n\to u\nabla u$ in $(L^2_{loc}(\RN))^N$ and $\nabla u_\n\to \nabla u$ a.e. on $\RN$.
Now we prove the conclusion in several steps.

\begin{itemize}[fullwidth,itemindent=0em]
\vskip 0.08in
\item[{\bf Step 1:}] We prove that $\nm{u_\n}_\iy\le C$ and $\nm{u}_\iy\le C$ for some positive constant $C$.

	We just prove the case $N\ge3$, the case $N=2$ can be obtained similarly.
	Set $T>2$, $r>0$ and 
	$$v_n=\begin{cases}T,\quad &u_n\ge T,\\ u_n,\quad &|u_n|\le T,\\-T,\quad &u_n\le -T. \end{cases}$$
Let $\phi=u_\n |v_n|^{2r}$, then $\phi\in\XR$. From $I_\n'(u_\n)+\la_\n u_\n=0$ and $\la_\n\ge0$,we obtain
\begin{align*}
	\int_\RN|u_\n|^{p-2}u_\n\phi
	&=\mu_\n\int_\RN|\nabla u_\n|^{\q-2}\nabla u_\n\cdot\nabla\phi +\int_\RN\nabla u_\n\cdot\nabla\phi\\
		&\quad~~   +2\int_\RN(u_\n\phi|\nabla u_\n|^2+|u_\n|^2\nabla u_\n\cdot\nabla\phi)+\la_\n\int_\RN u_\n\phi \\
	&\ge 2\int_\RN|u_\n|^2\nabla u_\n\cdot\nabla\phi \\
	&=2\int_\RN |u_\n|^2|\nabla u_\n|^2|v_n|^{2r}+|u_\n|^22r|v_n|^{2r-2}u_\n v_n\nabla u_\n\cdot\nabla v_n\\
	&=\f{1}{2}\int_\RN|v_n|^r|\nabla u_\n^2|^2+\f{4}{r}\int_\RN |u_\n^2\nabla |v_n|^r|^2\\
	&\ge\f{1}{r+4}\int_\RN|\nabla(u_\n^2|v_n|^2)|^2\ge \f{C}{(r+2)^2}\sbr{\int_\RN|u_\n^2|v_n|^2|^{2^*}}^{\f{2}{2^*}}.
\end{align*}
On the other hand, by the interpolation inequality, we have
\be\lab{tem102}
	\begin{aligned}
	\int_\RN|u_\n|^{p-2}u_\n\phi &=\int_\RN |u_\n|^p|v_n|^{2r}\\
		&\le\sbr{\int_\RN |u_\n|^{22^*}}^{\f{p-4}{22^*}} 
			\sbr{\int_\RN\sbr{|v_n|^r|u_\n|^2}^{\f{42^*}{22^*-p+4}}}^{\f{22^*-p+4}{22^*}}\\
		&\le C\sbr{\int_\RN\sbr{|v_n|^r|u_\n|^2}^{\f{42^*}{22^*-p+4}}}^{\f{22^*-p+4}{22^*}}.
	\end{aligned}
\ee
Combining these inequalities, one has
\be\lab{tem103}
	\sbr{\int_\RN|u_\n^2|v_n|^2|^{2^*}}^{\f{2}{2^*}}
		\le C(r+2)^2\sbr{\int_\RN\sbr{|v_n|^r|u_\n|^2}^{\f{42^*}{22^*-p+4}}}^{\f{22^*-p+4}{22^*}}.
\ee
Let $r_0:(r_0+2)q=22^*$ and $d=\f{2^*}{q}>1$ where $q=\f{42^*}{22^*-p+4}$. Taking $r=r_0$ in \eqref{tem103}, 
and letting $T\to+\iy$, we obtain
\be
	\nm{u_\n}_{(2+r_0)qd}\le \sbr{C(r_0+2)}^{\f{1}{r_0+2}}\nm{u_\n}_{(2+r_0)q}.
\ee
Set $2+r_{i+1}=(2+r_i)d$ for $i\in\N$. Then inductively, we have
\be\lab{tem104}
	\nm{u_\n}_{(2+r_0)qd^{i+1}}\le \prod_{k=0}^{i}\sbr{C(r_k+2)}^{\f{1}{r_k+2}}\nm{u_\n}_{(2+r_0)q}
			\le C_\iy \nm{u_\n}_{(2+r_0)q},
\ee
where $C_\iy$ is a positive constant. Taking $i\to\iy$ in \eqref{tem104}, we get
	$$\nm{u_\n}_\iy\le C \quad \text{and}\quad \nm{u}_\iy\le C.$$

\vskip 0.08in
\item[{\bf Step 2:}] We prove that $I'(u)+\la u=0$.

	Take $\phi=\psi e^{-u_\n}$ with $\psi\in\CR_0^\iy(\RN)$, $\psi\ge0$, we have 
\begin{align*}
	0&=\sbr{ I_\n'(u_\n)+\la_\n u_\n}[\phi]\\
	 &=\mu_n\int_\RN |\nabla u_\n|^{\q-2}\nabla u_\n\sbr{\nabla\psi e^{-u_\n}-\psi e^{-u_\n}\nabla u_\n}
	 				+\int_\RN \nabla u_\n\sbr{\nabla\psi e^{-u_\n}-\psi e^{-u_\n}\nabla u_\n}\\
	 &\quad~~  +2\int_\RN |u_\n|^2\nabla u_\n\sbr{\nabla\psi e^{-u_\n}-\psi e^{-u_\n}\nabla u_\n }
	 				+2\int_\RN u_\n\psi e^{-u_\n}|\nabla u_\n|^2\\
	 &\quad~~  +\la_\n \int_\RN u_\n\psi e^{-u_\n} -\int_\RN|u_\n|^{p-2}u_\n \psi e^{-u_\n}\\
	 &\le \mu_n\int_\RN |\nabla u_\n|^{\q-2}\nabla u_\n \nabla\psi e^{-u_\n} +\int_\RN(1+2u_\n^2)\nabla u_\n \nabla\psi e^{-u_\n}\\
	 &\quad~~  -\int_\RN (1+2u_\n^2-2u_\n)\psi e^{-u_\n}|\nabla u_\n|^2 
	 				+\la_\n \int_\RN u_\n\psi e^{-u_\n} -\int_\RN|u_\n|^{p-2}u_\n \psi e^{-u_\n}
\end{align*}
Since $\mu_n\to0^+$ and $\nm{u_\n}_\iy\le C$, \eqref{bdd} implies
	$$\n\int_\RN |\nabla u_\n|^{\q-2}\nabla u_\n \nabla\psi e^{-u_\n}\to0.$$
By the weak convergence of $u_\n$, the H\"older inequality and the Lebesgue's dominated convergence theorem we know that
	$$\int_\RN(1+2u_\n^2)\nabla u_\n \nabla\psi e^{-u_\n}\to\int_\RN(1+2u^2)\nabla u \nabla\psi e^{-u},$$
	$$\la_\n \int_\RN u_\n\psi e^{-u_\n}\to\la \int_\RN u\psi e^{-u},$$
and 
	$$\int_\RN|u_\n|^{p-2}u_\n \psi e^{-u_\n}\to\int_\RN|u|^{p-2}u \psi e^{-u}.$$
Moreover, by Fatou's lemma, there holds
	$$\liminf_{n\to\iy}\int_\RN (1+2u_\n^2-2u_\n)\psi e^{-u_\n}|\nabla u_\n|^2 \ge \int_\RN (1+2u^2-2u)\psi e^{-u}|\nabla u|^2.$$
Consequently, one has
\be\lab{tem20}
	\begin{aligned}
	0&\le\int_\RN \nabla u\sbr{\nabla\psi e^{-u}-\psi e^{-u}\nabla u} 
	 		+2\int_\RN |u|^2\nabla u\sbr{\nabla\psi e^{-u}-\psi e^{-u}\nabla u }\\
	 &\quad~~  +2\int_\RN u\psi e^{-u}|\nabla u|^2+\la_\n \int_\RN u\psi e^{-u} -\int_\RN|u|^{p-2}u \psi e^{-u}.
	\end{aligned}
\ee
For any $\vp\in\CR_0^\iy(\RN)$ with $\vp\ge0$. Choose a sequence of non-negative functions
 $\psi_n\in\CR_0^\iy(\RN)$ such that $\psi_n\to\vp e^{u}$ in $W^{1,2}(\RN)$, $\psi_n\to\vp e^{u}$ a.e. in $\RN$, and  that
$\psi_n$ is uniformly bounded in $L^\iy(\RN)$. Then we obtain from \eqref{tem20} that 
\be
	0\le\int_\RN \nabla u\cdot\nabla\vp+2\int_\RN (|u|^2\nabla u\cdot\nabla\vp+u\vp |\nabla u|^2)
				+\la \int_\RN u\vp -\int_\RN|u|^{p-2}u \vp.
\ee
Similarly by choosing $\phi=\psi e^{u_\n}$, we get an opposite inequality. Notice $\vp=\vp^+-\vp^-$ for any $\vp\in\CR_0^\iy(\RN)$, 
we get $I'(u)+\la u=0$.

\vskip 0.08in
\item[{\bf Step 3:}] We complete the proof.

	Similar as Lemma \ref{lempho}, we get from $I'(u)+\la u=0$ that
	$$Q(u):=Q_0(u)=0.$$
It follows that 
	$$Q_{\mu_n}(u_\n)+\ga_p\int_\RN|u_\n|^p\to Q(u)+\ga_p\int_\RN|u|^p.$$
Then using the weak lower semicontinuous property, there must be
\be\lab{tem21}
	\mu_n\int_\RN|\nabla u_\n|^\q\to0,\quad
	\int_\RN|\nabla u_\n|^2\to\int_\RN|\nabla u|^2,\quad
	\int_\RN|u_\n|^2|\nabla u_\n|^2\to\int_\RN|u|^2|\nabla u|^2.
\ee
That gives $I(u)=\lim_{n\to\iy}I_\mu(u_\n)=c$. Moreover, from \eqref{tem21}, we obtain
\be\lab{tem24}
	I_{\mu_n}'(u_\n)[u_\n]\to I'(u)[u].
\ee
Thus there holds $\la\nm{u_\n}_2^2\to\la\nm{u}_2^2$. 
So if $\la\neq0$, we have $\nm{u}_2^2=\lim_{n\to\iy}a_n$.
\end{itemize}
\ep

	Now we are able to end the proof of Theorem \ref{thm1} and \ref{thm2}. 
\bp[Proof of Theorem \ref{thm1} for the case $N\ge2$:]
	From Remark \ref{estofmin} and \ref{monotonicity}, we see that 
		$$d^*(a):=\lim_{\mu\to0^+}m_\mu(a)\in(0,+\iy).$$
By Theorem \ref{groundstate}, we can take $\mu_n\to0^+$, $I_\n'(u_\n)+\la_\n u_\n=0$ 
, $I_\n(u_\n)\to d^*(a)$ for $u_\n\in\SR_r(a_n)$ with $0<a_n\le a$ and $u_\n\ge0$. 
Then Lemma \ref{Liouville1} implies that $\la_\n>0$. 
Now Theorem \ref{converge} gives that there exist $v\neq0$, $v\ge0$, $v\in W^{1,2}_{rad}(\RN)\cap L^\iy(\RN)$ 
and $\la_0\in\R$ such that 
	$$\quad I'(v)+\la_0 v=0,\quad I(v)=d^*(a) \quad\text{and}\quad 0<\nm{v}_2^2\le a.$$
Thus by Lemma \ref{Liouville1}, there is $\la_0>0$. Since $\la_\n\to\la_0$,
 we may say that $\la_\n\neq0$ for $n$ large.
Then $a_n=a$ and $\nm{v }_2^2=a$. That is, $v$ is a nontrivial nonnegative solution of \eqref{massequ}. 
To consider the ground state normalized solution, we define 
	$$d(a):=\inf\lbr{I(v): v\in \tilde\SR(a), I|_{\tilde\SR(a)}'(v)=0, v\neq0}.$$
Then $d(a)\le I(v)=d^*(a)$. Futher, a similar approach to Lemma \ref{est} tells that $d(a)>0$.
We take a sequence $v_n\in\tilde\SR(a)$, $I|_{\tilde\SR(a)}'(v_n)=0$, $v_n\neq0$ and $v_n\ge0$ 
such that $I(v_n)\to d(a)$. We can show that (the proof is similar to that of Theorem \ref{converge}, so we omit it), 
up to a subsequence, there exist $u\neq0$, $u\ge0$, $u\in W^{1,2}_{rad}(\RN)\cap L^\iy(\RN)$ and $\la\in\R$ such that 
	$$\quad I'(u)+\la u=0\quad\text{and}\quad I(u)=d(a).$$
Again by Lemma \ref{Liouville1}, there is $\la\neq0$, and hence $\nm{u}_2^2=a$. 
That is, $u$ is a minimizer of $d(a)$. 
Finally, by \cite[Lemma 2.6]{Liu-Liu-Wang=JDE=2013}, $u$ is classical and strictly positive since $u\in L^\iy(\RN)$.
\ep

\bp[Proof of Theorem \ref{thm2}:]
	From Lemma \ref{genuslevel}, we see that 
		$$b_j(a)=\lim_{\mu\to0^+}c_\mu^j(a)\in(0,+\iy)\quad\text{and}\quad b_j(a)\to+\iy.$$
By Theorem \ref{infinite}, for each $j\in\N^+$ we can take $\nj\to0^+$, 
 $I_\nj'(u_\nj^j)+\la_\nj^j u_\nj^j=0$, $I_\nj(u_\nj^j)\to b_j(a)$ for $u_\nj\in\SR_r(a_n^j)$ with $0<a_n^j\le a$. 
And Lemma \ref{Liouville1} implies that $\la_\nj^j>0$. 
Now Theorem \ref{converge} gives that there exist $u^j\neq0$, $u^j\in W^{1,2}_{rad}(\RN)\cap L^\iy(\RN)$ 
and $\la^j\in\R$ such that 
	$$\quad I'(u^j)+\la^j u^j=0,\quad I(u^j)=b_j(a) \quad\text{and}\quad 0<\nm{u^j}_2^2\le a.$$
Thus by Lemma \ref{Liouville1}, there is $\la^j>0$. Going back since $\la_\nj^j\to\la^j$,
 we may say that $\la_\nj^j\neq0$ for $n$ large. 
Then $a_n^j=a$ and $\nm{u^j }_2^2=a$. That is $\lbr{u^j:j\in\N^+}$ is a sequence of normalized solutions of \eqref{massequ}. 
Moreover, $I(u^j)=b_j\to+\iy$. 
\ep

\vskip0.23in
\section{The mass critical case \texorpdfstring{$p=4+\f{4}{N}$}{} }

	In this section we denote $p_*=4+\f{4}{N}$  and assume that $p=p_*$. 
We still consider $I_\mu$, but on an open subset of $\XR$. Let
\be
	\OR:=\lbr{u\in\XR:\int_\RN u^2|\nabla u|^2<\f{N}{4(N+1)}\int_\RN |u|^{p_*} },
\ee
and for simplicity, we still denote
	$$\SR(a):=\lbr{u\in\OR:\int_\RN u^2=a},$$
	$$\QR_\mu(a):=\lbr{u\in\SR(a):Q_\mu(u)=0},$$
	$$\SR_r(a):=\SR(a)\cap \XR_r,\quad \QR_\mu^r(a):=\QR_\mu(a)\cap\XR_r.$$
We have 
\bl\lab{nonempty}
	When $a>a^*$, $\SR(a)$ is nonempty.
\el
\bp
	Let $u=Q_{p_*}^{\f{1}{2}}$, then from \eqref{SGNine}, we have
\be\lab{tem201}
	\int_\RN |u|^{p_*}=\f{4(N+1)}{N}\int_\RN u^2|\nabla u|^2.
\ee
Let $w_a=\sbr{\f{a}{a_*}}^{\f{1}{2}}u$, then $\nt{w_a}^2=a$ and \eqref{tem201} implies that
\be
	\int_\RN w_a^2|\nabla w_a|^2=\f{N}{4(N+1)}\sbr{\f{a}{a_*}}^{-\f{2}{N}}\int_\RN |w_a|^{p_*}
		<\f{N}{4(N+1)}\int_\RN |w_a|^{p_*},
\ee
that is $w_a\in\SR(a)$.
\ep

	So from now on, we assume $a>a^*$. 
Then noting that when $p=p_*$, there is $p_*\ga_{p_*}>\q+\q\ga_\q$ and $p_*\ga_{p_*}=2+N$, we still have

\bl
	Let $0<\mu\le1$, then $\QR_\mu(a)$ is a $\CR^1$-submanifold of codimension 1 in $\SR(a)$, hence a $\CR^1$-submanifold of 
codimension 2 in $\XR$.
\el

\bl\lab{structure1}
	For any $0<\mu\le1$ and any $u\in\OR\setminus\{0\}$, the following statements hold.
	\begin{itemize}[fullwidth,itemindent=0em]
	\item[(1)]	There exists a unique number $s_\mu(u)\in\R$ such that $Q_\mu(s_\mu(u)\s u)=0$.
	\item[(2)]	$I_\mu(s\s u)$ is strictly increasing in $s\in(-\iy,s_\mu(u))$ and is strictly decreasing 
	in $s\in(s_\mu(u),+\iy)$, and
		 $$\lim_{s\to-\iy}I_\mu(s\s u)=0^+,\quad \lim_{s\to+\iy}I_\mu(s\s u)=-\iy,\quad I_\mu(s_\mu(u)\s u)>0.$$
	\item[(3)]	$s_\mu(u)<0$ if and only if $Q_\mu(u)<0$.
	\item[(4)]  The map $u\in\XR\setminus\{0\} \mapsto s_\mu(u)\in\R$ is of class $\CR^1$.
	\item[(5)]  $s_\mu(u)$ is an even function with respect to $u\in \XR\setminus\{0\}$.
	\end{itemize}
\el

	Similarly to Lemma \ref{est}, there also holds
\bl\lab{est1}
	The following statements hold.
	\begin{itemize}[fullwidth,itemindent=0em]
	\item[(1)]	$\DR(a):=\inf_{0<\mu\le1,u\in\QR_\mu(a)}\int_\RN|u|^2\gu{u}2>0$ is independent of $\mu$.
	\item[(2)]	If $\sup_{n\ge1}I_\mu(u_n)<+\iy$ for $u_n\in\QR_\mu(a)$, then
			$$\sup_{n\ge1} \max \lbr{ \mu\int_\RN\gu{u}\q, \int_\RN|u|^2\gu{u}2, \int_\RN\gu{u}2}<+\iy.$$
	\end{itemize}
\el
\bp
	The proof is different from the one of Lemma \ref{est}.
\begin{itemize}[fullwidth,itemindent=0em]
\item[(1)] For any $u\in\QR_\mu(a)$, using the equality $Q_\mu(u)=0$ and \eqref{GNine} we obtain
\be\lab{control}
	\int_\RN |\nabla u|^2\le(N+2)\mbr{\sbr{\f{a}{a_*}}^{\f{2}{N}}-1}\int_\RN u^2|\nabla u|^2.
\ee
On the one hand, when $N\le3$ , there holds $p_*<2^*$. 
Therefore the calssical Gagliardo-Nirenberg inequality (\cite{GNinequality}) tells that 
\be
	\int_\RN|\nabla u|^2\le\ga_{p_*}\int_\RN |u|^{p_*}\le C(N) a^{1+\f{2}{N}-\f{N}{2}}\sbr{\int_\RN |\nabla u|^2}^{\f{N+2}{2}},
\ee
follows which there is $\int_\RN|\nabla u|^2\ge\f{C(N)}{a^{\f{4}{N^2}+\f{2}{N}-1}}$.
Combining with \eqref{control}, one obtain 
	$$\inf_{0<\mu\le1,u\in\QR_\mu(a)}\int_\RN|u|^2\gu{u}2>0.$$
On the other hand, when $N\ge4$, there is $p_*>2^*$. But using interpolation inequality and Young inequality we have
\be
	\begin{aligned}
	(N+2)\int_\RN u^2|\nabla u|^2+\int_\RN|\nabla u|^2 &\le \ga_{p_*}\int_\RN |u|^{p_*}
				\le\sbr{ \int_\RN |u|^{2^*}}^{\f{22^*-p_*}{2^*}}\sbr{\int_\RN|u|^{22^*}}^{\f{p_*-2^*}{2^*}}\\
		&\le C(N)\sbr{ \int_\RN |\nabla u|^2}^{\f{22^*-p_*}{2}}\sbr{\int_\RN u^2|\nabla u|^2}^{\f{p_*-2^*}{2}}\\
		&\le (N+2)\int_\RN u^2|\nabla u|^2 +C(N)\sbr{ \int_\RN |\nabla u|^2} ^{\f{22^*-p_*}{2^*+2-p_*}},
	\end{aligned}
\ee
which gives that $\int_\RN|\nabla u|^2\ge C(N)$ and again
	$$\inf_{0<\mu\le1,u\in\QR_\mu(a)}\int_\RN|u|^2\gu{u}2>0.$$

\vskip 0.08in
\item[(2)]	Since $p_*\ga_{p_*}=2+N$, we see from \eqref{minus} that
	$$\sup_{n\ge1} \max \lbr{ \mu\int_\RN\gu{u}\q, \int_\RN\gu{u}2}<+\iy.$$
Moreover, $Q_\mu(u_n)=0$ implies that
\be
	\begin{aligned}
	C&\ge \mu\int_\RN\gu{u}\q+\int_\RN\gu{u}2\\
	 &= \ga_{p_*}\int_\RN|u_n|^{p_*}-(2+N)\int_\RN u_n^2|\nabla u_n|^2\\
	 &> (2+N)\sbr{\f{N}{4(N+1)}-1}\int_\RN u_n^2|\nabla u_n|^2,
	\end{aligned}
\ee
which finishes  the proof. 
\end{itemize}
\ep

First, we will consider a minimizition problem
\be\lab{min1}
	m_\mu(a):=\inf_{u\in\QR_\mu(a)} I_\mu(u).
\ee
\br
It is easy to see from Lemma \ref{est1} and \eqref{minus} that 
\be
	\inf_{0\le\mu\le1} m_\mu(a)\ge\f{N}{2(2+N)}\inf_{0\le\mu\le1,u\in\QR_\mu(a)}\int_\RN\gu{u}2>0.
\ee
On the other hand, to use  the convergence Theorem \ref{converge}, we need to give an uniform upper bound of $m_\mu(a)$. 
Indeed for any fixed $a>a^*$, recalling the function $w_a=\sbr{\f{a}{a_*}}^{\f{1}{2}}Q_{p_*}^{\f{1}{2}}\in\SR(a)$ 
in Lemma \ref{nonempty}, and let $s_\mu:=s_\mu(w_a)$, then from $Q_\mu(s_\mu\s w_a)=0$ we obtain
\be
	\begin{aligned}
	&\quad~~ (1+\ga_\q)\mu e^{-(2+N-\q-\q\ga_\q)s_\mu}\sbr{\f{a}{a_*}}^{\f{\q}{2}}\int_\RN|\nabla Q_{p_*}^{\f{1}{2}}|^\q
			+e^{-Ns_\mu}\sbr{\f{a}{a_*}}\int_\RN|\nabla Q_{p_*}^{\f{1}{2}}|^2\\
	&=(1+\ga_\q)\mu e^{-(2+N-\q-\q\ga_\q)s_\mu}\int_\RN|\nabla w_a|^\q+e^{-Ns_\mu}\int_\RN|\nabla w_a|^2\\
	&=\ga_{p_*}\int_\RN|w_a|^{p_*}-(2+N)\int_\RN|w_a|^2|\nabla w_a|^2\\
	&=\f{N(2+N)}{4(N+1)}\sbr{ 1-\sbr{\f{a}{a_*}}^{-\f{2}{N}} }\sbr{\f{a}{a_*}}^{2+\f{2}{N}}\nm{Q_{p_*}^{\f{1}{2}}}_1>0,
	\end{aligned} 
\ee
it follows that   $\sup_{0\le\mu\le1}s_\mu<+\iy$. Therefore,
\be
	\begin{aligned}
	\sup_{0\le\mu\le1} m_\mu(a)&\le \sup_{0\le\mu\le1} I_\mu(s_\mu\s w_a)
				=\sup_{0\le\mu\le1} I_\mu(s_\mu\s w_a)-Q_\mu(s_\mu\s w_a)\\
		&=\sup_{0\le\mu\le1} \f{2+N-\q-\q\ga_\q}{\q (2+N)}\mu e^{\q(1+\ga_\q)s_\mu}\int_\RN|\nabla Q_{p_*}^{\f{1}{2}}|^\q
			+\f{N}{2(2+N)}e^{2s_\mu}\int_\RN|\nabla Q_{p_*}^{\f{1}{2}}|^2\\
		&< +\iy.
	\end{aligned}
\ee
\er
Now we construct a special Palais-Smale sequence of $I_\mu|_{\SR(a)}$ at level $m_\mu(a)$. 
But different from  the Section \ref{groundstatesolution}, in mass-critical case there is no result  as Lemma \ref{tem4},
and hence there is no the mountain-pass type result  as  Lemma \ref{mountain-pass}. 
So we will not consider $I_\mu$ directly. 
Instead we study the auxilary functional $K_\mu(u)$ defined by \eqref{auxilary1}, 
and we point out that our approach is inspired by \cite{BS=CVPDE=2019} 
(see also \cite{Jeanjean=TransAMS=2019}). Similar to \cite[Lemma 3.7]{BS=CVPDE=2019}, we have
\bl\lab{welldefined}
	Let a sequence $u_n\in\SR(a)$ with $u_n\to u$ in $\XR$ as $n\to\iy$. 
	Then if $u\in\pa\OR$, we have $K_\mu(u_n)\to\iy$ as $n\to\iy$.
\el
\bp
	If $u_n\to u$ in $\XR$, then there are
		$$\int_\RN|\nabla u_n|^\q \to\int_\RN\gu{u}\q>0,\quad \int_\RN|\nabla u_n|^2 \to \int_\RN\gu{u}2>0,$$
		$$\int_\RN|u_n|^2|\nabla u_n|^2 \to\int_\RN|u|^2\gu{u}2>0,\quad \int_\RN|u_n|^{p_*}\to\int_\RN|u|^{p_*}>0.$$
Let $s_n=s_\mu(u_n)$. Since $Q_\mu(s_n\s u_n)=0$, we obtain
\be
	\begin{aligned}
	&\quad~~ (1+\ga_\q)\mu e^{-(2+N-\q-\q\ga_\q)s_n}\int_\RN|\nabla u_n|^\q+e^{-Ns_n}\int_\RN|\nabla u_n|^2\\
	&=\ga_{p_*}\int_\RN|u_n|^{p_*}-(2+N)\int_\RN|u_n|^2|\nabla u_n|^2\\
	&\to \ga_{p_*}\int_\RN|u|^{p_*}-(2+N)\int_\RN|u|^2|\nabla u|^2=0,
	\end{aligned}
\ee
where the last equality comes from $u\in\pa\OR$. It follows that $s_n\to+\iy$. So
\be
	\begin{aligned}
	K_\mu(u_n)&=I_\mu(s_n\s u_n)=I_\mu(s_n\s u_n)-Q_\mu(s_n\s u_n)\\
		&=\f{2+N-\q-\q\ga_\q}{\q (2+N)}\mu e^{\q(1+\ga_\q)s_n}\int_\RN|\nabla u_n|^\q
			+\f{N}{2(2+N)}e^{2s_n}\int_\RN|\nabla u_n|^2\\
		&\to+\iy.
	\end{aligned}
\ee
\ep

	Recalling the Definition A in Section \ref{groundstatesolution}, we give directly the following results without proof, 
since the proof is very similar to the one of \cite[Proposition 3.9]{BS=CVPDE=2019} 
(see also \cite{Jeanjean=TransAMS=2019}).
\bl\lab{minmax1}
	Let $\GR$ be a homotopy  stable family of compact subsets of $Y=\SR_r(a)$ with boundary $B=\emptyset$, 
	and set
		$$d:=\inf_{A\in\GR}\max_{u\in A}K_\mu(u).$$
	If $d>0$,  then there exists a sequence $u_n\in \SR_r(a)$ such that as $n\to\iy$,
		$$I_\mu(u_n)\to d,\quad I_\mu|_{\SR(a)}'(u_n)\to0, \quad Q_\mu(u_n)=0.$$
	Moreover, if one can find a minimizing sequence $A_n$ for $d$ with the property that $u\ge0$ a.e. for any $u\in A_n$,
	then one can find the sequence $u_n$ satisfying the additional condition
		$$u_n^-\to0,\quad \text{a.e. in}~\RN.$$
\el
\br
	As pointed out in \cite{BS=CVPDE=2019}, the set $\OR$ is neither complete, nor connected, and hence 
	in principle the assumptions of the minimax theorem (such as \cite[Theorem 3.2]{Ghoussoub=1993}) are not satisfied. 
	However, the connectedness assumption can be avoided considering the restriction of $K_\mu$ on the connected component of $\OR$
	(if $B\neq\emptyset$, we need to assume that $B$ is contained in a connected component of $\QR_\mu(a)$).
	Regarding the completeness, what is really used in the deformation lemma \cite[Lemma 3.7]{Ghoussoub=1993} is that 
	the sublevel sets $K_\mu^c:=\lbr{u\in\SR(a): K_\mu(u)\le c}$ are complete for every $c\in\R$. 
	This follows by Lemma \ref{welldefined}. Hence the minmax theorem \cite[Theorem 3.2]{Ghoussoub=1993} can be used to
	obtain the Palais-Smale sequence.  The rest of the process is similar to Lemma \ref{minmax}. 
\er

\bl
	For any fixed $\mu\in(0,1]$, there exists a sequence $u_n\in\SR_r(a)$ such that
		$$I_\mu(u_n)\to m_\mu(a),\quad I_\mu|_{\SR(a)}'(u_n)\to0, 
		\quad Q_\mu(u_n)=0\quad\text{and}\quad u_n^-\to0~\text{a.e. in}~\RN.$$
\el
\bp
	We use Lemma \ref{minmax1} by taking the set $\GR$ of all singletons belonging to $\SR_r(a)$. 
It is clearly a homotopy stable family of compact subsets of $\SR_r(a)$ with boundary $B=\emptyset$. 
Observe that 
	$$\al_\mu(a)=\inf_{A\in\GR}\max_{u\in A}K_\mu(u)=\inf_{u\in\SR_r(a)}\max_{s\in\R}I_\mu(s\s u).$$
We claim that 
	$$\al_\mu(a)=m_\mu(a).$$
Indeed, on the one hand, for any $u\in\SR_r(a)$ there exists a $s_\mu(u)$ such that $s_\mu(u)\s u\in\QR_\mu(a)$ 
and $I_\mu(s_\mu(u)\s u)=\max_{s\in\R}I_\mu(s\s u)$. This implies that 
	$$\al_\mu(a)=\inf_{u\in\SR_r(a)}\max_{s\in\R}I_\mu(s\s u)\ge \inf_{u\in\QR_\mu(a)}I_\mu(u)=m_\mu(a).$$
On the other hand, for any $u\in\QR_\mu^r(a)$, $I_\mu(u)=\max_{s\in\R}I_\mu(s\s u)$, so
	$$m_\mu^r(a):=\inf_{u\in\QR_\mu^r(a)}I_\mu(u)\ge \inf_{u\in\SR_r(a)}\max_{s\in\R}I_\mu(s\s u)=\al_\mu(a).$$
Finally the inequality $m_\mu(a)\ge m_\mu^r(a)$ can be obtained easily by the symmetric decreasing rearrangement.
Thus, the conclusion follows directly from Lemma \ref{minmax1}.
\ep

	Then as in Section \ref{groundstatesolution}, we have 
\bt\lab{groundstate1}
	Let $p=p_*$. For any fixed $\mu\in(0,1]$, there exists a $u_\mu\in\XR_r\setminus\{0\}$ and a $\la_\mu\in\R$ such that 
\begin{align*}
	&I_\mu'(u_\mu)+\la_\mu u_\mu=0,\\
	&I_\mu(u_\mu)=m_\mu(a),\quad Q_\mu(u_\mu)=0,\\
	&0<\nm{u_\mu}_2^2\le a,\quad u_\mu\ge0.
\end{align*}
Moreover, if $\la_\mu\neq0$, we have that $\nm{u_\mu}_2^2=a$, i.e.,  $m_\mu(a)$ is achieved, and $u_\mu$ is 
a ground state critical point of $I_\mu|_{\SR(a)}$.
\et

\bp[Proof of Theorem \ref{thm3}: ]
	The proof is exactly the same as the one of Theorem \ref{thm1}, so we omit the details.
\ep
\br
	We are not able to obtain multiple solutions as in Section \ref{infinitelymanysolutions}. 
Indeed, if we consider  an open subset $\OR$ and follow the strategy in Section \ref{infinitelymanysolutions}, 
we need to prove a result like Lemma \ref{genuslevel}. 
However for any   finite dimensional subspace $W_j$ of $\XR$, using the equivalence of norms in finite dimensional spaces, 
we can only obtain that for any $j>0$, there exists a $a(j)>0$ large enough such that 
	$$\lbr{u\in W_j:\nt{u}^2=a}\subset \OR\quad \text{when}\quad a>a(j),$$
which is necessary to prove the nonemptyness of the sets of type $\AR_j$.
And another difficulty is that as $\mu\to0^+$, we are unable to distinguish the energy 
	$$b_j(a):=\lim_{\mu\to0^+}c_\mu^j(a)\quad \text{and}\quad b_k(a):=\lim_{\mu\to0^+}c_\mu^k(a),$$
for $j\neq k$. As a result, we can not distinguish the solutions related to $b_j(a)$ and $b_k(a)$.
\er

	Recalling Proposition \ref{prop1}, we prove the concentration theorem.
\bp[Proof of the Theorem \ref{thm4}: ]
	Let $u_n$ be a radially symmetric positive solution of \eqref{massequ} for $a=a_n$ with $a_n>a_*$ and $a_n\to a_*$.
From Lemma \ref{est1}, we see that 
\be
	\int_\RN u_n^2|\nabla u_n|^2\ge \f{C}{\sbr{\f{a_n}{a_*}}^{\f{2}{N}}-1 }\to+\iy,
\ee
\be\lab{tem51}
	\f{\int_\RN|\nabla u_n|^2}{\int_\RN u_n^2|\nabla u_n|^2}\le C\sbr{\sbr{\f{a_n}{a_*}}^{\f{2}{N}}-1}\to0.
\ee
Since $Q_\mu(u_n)=0$, we know that
\be\lab{tem52}
	\f{\int_\RN|u_n|^{p_*}}{\int_\RN u_n^2|\nabla u_n|^2}\to \f{4(N+1)}{N}.
\ee
Let $v_n(x):=\e_n^{\f{N}{2}}u_n(\e_n x)$ with
	$$\e_n=\sbr{ \int_\RN u_n^2|\nabla u_n|^2 }^{-\f{1}{2+N}}\to0^+.$$
Direct calculations show that $\nt{v_n}^2=a_n\to a_*$, $\int_\RN v_n^2|\nabla v_n|^2=1$, 
$\nm{v_n}_{p_*}^{p_*}\to \f{4(N+1)}{N}$ and $\e_n^N\nt{\nabla v_n}^2\to0$.
Then $v_n^2$ is bounded in $\ER^{p_*}$. Moreover, using \cite[Lemma I.1]{Lions2}, we deduce that there exist $\delta>0$ 
and a sequence $y_n\in\RN$ such that for some $R>0$,
	$$\int_{B_R(y_n)} v_n^2\ge\delta.$$
Thus there exists a nonnegative radially symmetric function $v\neq0$ with $v^2\in\ER^{p_*}\cap L^2(\RN)$ such that
	$$v_n^2(\cdot +y_n)\rh v^2    \quad\text{in}~\ER^{p_*},$$
	$$v_n(\cdot +y_n)\rh v        \quad\text{in}~L^2(\RN),$$
	$$v_n^2(\cdot +y_n)\to v^2    \quad\text{in}~L^q(\RN)~\text{for}~1<q<2^*,$$
	$$v_n(\cdot +y_n)\to v        \quad\text{a.e. in}~\RN.$$
Since $u_n$ solves 
	$$-\Delta u_n-u_n\Delta u_n^2+\la_n u_n=u_n^{p_*-1},$$
where the Lagrange multiplier is given by 
	$$\la_n=\f{1}{a_n}\sbr{ \int_\RN|u_n|^{p_*}-\int_\RN|\nabla u_n|^2-\int_\RN u_n^2|\nabla u_n|^2 },$$
and $v_n$ satisfies
	$$-\e_n^N\Delta v_n-v_n\Delta v_n^2+\e_n^{2+N}\la_n v_n=v_n^{p_*-1}.$$
Combining \eqref{tem51} and \eqref{tem52}, we deduce that $\e_n^{2+N}\la_n\to \f{4}{Na^*}$.
Then a similar approach as Lemma \ref{A2} tells that 
\be\lab{tem53}
	-v\Delta v^2+\e_n^{2+N}\la_n v=v^{p_*-1}.
\ee
Now setting
\be
	\begin{aligned}
	w_n(x):&=\sbr{ \f{Na^*}{4} }^{\f{N}{2+N}}  v_n^2\sbr{ \sbr{ \f{Na^*}{4} }^{\f{1}{2+N}}x+y_n }\\
		&=\mbr{ \sbr{ \f{Na^*}{4} }^{\f{1}{2+N}}\e_n}^N  u_n^2\sbr{ \sbr{ \f{Na^*}{4} }^{\f{1}{2+N}}\e_nx + \e_ny_n },
	\end{aligned}
\ee
\be
	w(x):=\sbr{ \f{Na^*}{4} }^{\f{N}{2+N}}  v^2\sbr{ \sbr{ \f{Na^*}{4} }^{\f{1}{2+N}}x },
\ee
it is easily seen that $w_n\rh w$ in $\ER^{p_*}$ and $\nm{w_n}_1=\nt{v_n}^2=a_n$. Moreover, it follows from \eqref{tem53} 
that $w$ is a solution of \eqref{GNequ1}. Thus $w=Q_{p_*}$, and hence $\nm{w}_1=\nt{v}^2=a_*$. 
So we have $v_n\to v$ in $L^2(\RN)$, which finishes  the proof. 
\ep

\vskip0.26in
\appendix
\renewcommand{\appendixname}{Appendix \Alph{section}}
\section{Appendix}

\bl\lab{A1}
	In the setting of Section \ref{sec1}, $V(u)\in\CR^1(\XR)$. 
\el
\bp
	The proof is elementry. When $N=2$, since $W^{1,\q}(\R^2)\hra \CR^{0,\al}(\R^2)$, it is easily to check
that $V(u)\in\CR^1(\XR)$. Now we set $N\ge3$. For any $u,\phi\in\XR$, 
\be\lab{temA1}
	\f{V(u+t\phi)-V(u)}{t}=At+Bt^2+Ct^3+2\int_\RN u\phi|\nabla u|^2+u^2\nabla u\cdot\nabla \phi,
\ee
where
	$$A=\int_\RN u^2|\nabla\phi|^2+\phi^2|\nabla u|^2+4u\phi\nabla u\cdot\nabla\phi,$$
	$$B=\int_\RN\phi^2\nabla u\cdot\nabla\phi+u\phi|\nabla\phi|^2 \quad\text{and}\quad
		C=\int_\RN\phi^2|\nabla\phi|^2.$$
We need to prove $A,B,C$ are finite numbers. Indeed, since $\f{4N}{N+2}<\q<\f{4N+4}{N+2}<4$, there is 
$\q<\f{2\q}{\q-2}< \f{\q N}{N-\q}$ and hence
\be\lab{temA2}
	\int_\RN u^2|\nabla\phi|^2\le \sbr{\int_\RN |u|^{\f{2\q}{\q-2}}}^{\f{\q-2}{\q}}\sbr{\int_\RN |\nabla\phi|^\q}^{\f{2}{\q}}
		\le C\nm{u}_{W^{1,\q}(\RN)}^{\f{2}{\q}} \nm{\phi}_{W^{1,\q}(\RN)}^{\f{2}{\q}}<\iy,
\ee
and we can handle other terms in a similar way, so $A,B,C$ are finite numbers. 
Now by letting $t\to0$ in \eqref{temA1}, we immediately get the Fr\`echet deravetive is
	$$DV(u)[\phi]=2\int_\RN u\phi|\nabla u|^2+u^2\nabla u\cdot\nabla \phi.$$
Then in a similarly way in \eqref{temA2}, one can prove that $DV(u)$ is continuous for $u\in\XR$, 
so  $V(u)\in\CR^1(\XR)$ and $V'(u)=DV(u)$.
\ep

\bl\lab{A2}
	Assume that $I_\mu'(u_n)+\la u_n\to0$ for some $\la\in\R$ with $u_n\in\XR$, and that $u_n\rh u$ in $\XR$. Then up to a subsequence, 
\begin{itemize}[fullwidth,itemindent=0em]
\item[(1)]	$u_n\to u$ in $\XR_{loc}:=W^{1,\q}_{loc}(\RN)\cap W^{1,2}_{loc}(\RN)$,
\item[(2)]	$u_n\nabla u_n\to u\nabla u$ in $(L^2_{loc}(\RN))^N$,
\item[(3)]	$I_\mu'(u)+\la u=0$.
\end{itemize}
\el
\bp
	The proof is inspired by \cite[Lemma 14.3]{Kuzin-Pohozaev=1997}. 
Since $u_\n\rh u$ in $\XR$, we have $\sup_{n} \nm{u_n}_\XR\le C_0$. 
For any $R>1$, we set $\phi\in\CR_0^\iy(\RN)$ satisfying
	$$0\le\phi\le1,\quad \phi(x)=\begin{cases} 1,\quad &|x|\le R,\\0,\quad &|x|\ge 2R,\end{cases} 
		\quad\text{and}\quad  |\nabla\phi|\le2.$$
Then for any $n,m\in\N$,
\be\lab{temA3}
	\begin{aligned}
	o(1)_n+o(1)_m&=(I_\mu'(u_n)+\la u_n)[(u_n-u_m)\phi]-(I_\mu'(u_m)+\la u_m)[(u_n-u_m)\phi]\\
		&=\mu\int_\RN (|\nabla u_n|^{\q-2}\nabla u_n-|\nabla u_m|^{\q-2}\nabla u_m)\cdot\nabla\sbr{(u_n-u_m)\phi}\\
			&\quad~~ +\int_\RN (\nabla u_n-\nabla u_m)\cdot\nabla\sbr{(u_n-u_m)\phi}\\
			&\quad~~ +2\int_\RN (u_n|\nabla u_n|^2-u_m|\nabla u_m|^2)(u_n-u_m)\phi\\
			&\quad~~ +2\int_\RN (u_n^2\nabla u_n-u_m^2\nabla u_m)\cdot\nabla\sbr{(u_n-u_m)\phi}\\
			&\quad~~ -\int_\RN (|u_n|^{p-2}u_n-|u_m|^{p-2}u_m)(u_n-u_m)\phi\\
		&=:K_1+K_2+K_3+K_4+K_5.
	\end{aligned}
\ee
Next we estimate $K_i$ for $i=1,2,3,4,5$ one by one.
\begin{align*}
	K_1&=\mu\int_{B_R}(|\nabla u_n|^{\q-2}\nabla u_n-|\nabla u_m|^{\q-2}\nabla u_m)\cdot\nabla(u_n-u_m)\\
			&\quad~~ +\mu\int_{B_{2R}\setminus B_R}(|\nabla u_n|^{\q-2}\nabla u_n-|\nabla u_m|^{\q-2}\nabla u_m)
						\cdot\nabla(u_n-u_m)\phi\\
			&\quad~~ +\mu\int_{B_{2R}\setminus B_R}(|\nabla u_n|^{\q-2}\nabla u_n-|\nabla u_m|^{\q-2}\nabla u_m)
						\cdot\nabla\phi (u_n-u_m)\\
		&\ge C\mu\int_{B_R}|\nabla u_n-\nabla u_m|^\q+C\mu\int_{B_{2R}\setminus B_R}|\nabla u_n-\nabla u_m|^\q\phi\\
			&\quad~~ -C\sbr{\nm{u_n}_\q^{\q-1}+\nm{u_m}_\q^{\q-1}} \nm{u_n-u_m}_{L^\q(B_{2R})}\\
		&\ge C\mu\nm{\nabla u_n-\nabla u_m}_{L^\q(B_R)}^\q -C\nm{u_n-u_m}_{L^\q(B_{2R})},
\end{align*}
and similarly
	$$K_2\ge C\nm{\nabla u_n-\nabla u_m}_{L^2(B_R)}^2 -C\nm{u_n-u_m}_{L^2(B_{2R})},$$
	$$K_3\ge -C\nm{u_n-u_m}_{L^\q(B_{2R})},$$
	$$K_4\ge 2\nm{u_n\nabla u_n-u_m\nabla u_m}_{L^2(B_R)}^2-C\nm{u_n-u_m}_{L^\q(B_{2R})},$$
	$$K_5\ge -C\nm{u_n-u_m}_{L^p(B_{2R})}.$$
Substituting these estimates into \eqref{temA3}, we obtain
\be
	\begin{aligned}
	&\quad~ \mu\nm{\nabla u_n-\nabla u_m}_{L^\q(B_R)}^\q+\nm{\nabla u_n-\nabla u_m}_{L^2(B_R)}^2
			+\nm{u_n\nabla u_n-u_m\nabla u_m}_{L^2(B_R)}^2\\
	&\le C\nm{u_n-u_m}_{L^\q(B_{2R})}+C\nm{u_n-u_m}_{L^2(B_{2R})}
			+C\nm{u_n-u_m}_{L^\q(B_{2R})}+o(1)_n+o(1)_m\\
	&\to0,\quad \text{as}~n\to\iy,~m\to\iy,
	\end{aligned}
\ee
where in the last estimate we use the compact embedding theorem in bounded domains.
Thus for any $R>1$, $u_n$ is a Cauchy sequence in $W^{1,\q}(B_R)\cap W^{1,2}(B_R)$,  and $u_n\nabla u_n$ is also a 
Cauchy sequence in $\sbr{L^2(B_R)}^N$. So up to a subsequence $u_n\to u$ in $\XR_{loc}$ and $u_n\nabla u_n\to u\nabla u$ 
in $\sbr{L^2_{loc}(\RN)}^N$. 
Finally, we need to prove that for any $\vp\in\XR$, there holds $(I_\mu'(u)+\la u)[\vp]=0$. 
Since $u_n\nabla u_n\to u\nabla u$ a.e. in $\RN$ and $u_n$ is bounded in $\XR$, we obtain that
	$$|\nabla u_n|^{\q-2}\nabla u_n\rh |\nabla u|^{\q-2}\nabla u\quad \text{in}~L^{\f{\q}{\q-1}}(\RN),$$
	$$u_n|\nabla u_n|^2\rh u|\nabla u|^2\quad \text{in}~L^{\f{4}{3}}(\RN),$$
	$$u_n^2 \nabla u_n\rh u^2\nabla u \quad \text{in}~\sbr{L^{\f{4}{3}}(\RN)}^N,$$
it follows that 
	$$(I_\mu'(u)+\la u)[\vp]=\lim_{n\to\iy}(I_\mu'(u_n)+\la u_n)[\vp]=0.$$
\ep



\vskip0.26in









\vskip0.26in
\begin{thebibliography}{10}

 
 \bibitem{Agueh=2008}
M. Agueh. 
Sharp Gagliardo-Nirenberg inequalities via $p$-Laplacian type equations.
{\it Nonliear Differential Equations Appl.}, {\bf 15}(2008), 457-472.

 \bibitem{BJS=JMPA=2016}
T.~Bartsch, L.~Jeanjean and N.~Soave.
\newblock Normalized solutions for a system of coupled cubic Schr\"odinger equations on $\R^3$.
{\it J. Math. Pures Appl.}, {\bf (9)106}(2016), no.4, 583-614.

\bibitem{BartschJeanjean=2018}
T.~Bartsch and L.~Jeanjean.
\newblock Normalized solutions for nonlinear Schr\"odinger systems.
{\it Proc. Roy. Soc. Edinburgh Sect. A}, {\bf 148}(2018), no.2, 225-242.



\bibitem{Bartsch-Valeriola=2013}
T. Bartsch and S. de Valeriola. 
Normalized solutions of nonlinear Schr\"odinger equations. 
{\it Arch. Math.}, {\bf 100}(2013), 75-83.

\bibitem{BartschSoave=2017}
T.~Bartsch and N.~Soave. 
A natural constraint approach to normalized solutions of nonlinear Schr\"odinger equations and systems. 
{\it J. Funct. Anal.}, {\bf 272}(2017), no.12, 4998-5037.

 
 
 
\bibitem{BS=CVPDE=2019}
T.~Bartsch and N.~Soave.
\newblock Multiple normalized solutions for a competing system of Schr\"odinger equations.
{\it Calc. Var. Partial Differential Equations}, {\bf 58}(2019), no.1, Art.22, 24 pp.


\bibitem{Zhong-Zou-Bartsch=2020}
T. Bartsch, X. X. Zhong and W. M. Zou. 
Normalized solutions for a coupled Schr\"odinger system. 
{\it Math. Ann.} (2020). https://doi.org/10.1007/s00208-020-02000-w.
 

 

\bibitem{bg-1}
G. Bass and N. Nasonov. 
Nonlinear electromagnetic-spin waves. 
{\it Phys. Rep.}, {\bf 189}(1990), 165-223.
 
 
 
 
 
\bibitem{Bellazzini-Jeanjean-Luo=2013}
J. Bellazzini, L. Jeanjean and T. Luo. 
Existence and instability of standing waves with prescribed norm for a class of Schr\"odinger–Poisson equations.  
{\it Proc. Lond. Math. Soc.}, {\bf 107}(2013), 303-339.


\bibitem{Berestycki-Lions=1983-1}
H. Berestycki and P. L. Lions. 
Nonlinear scalar field equations I: Existence of a ground state. 
{\it Arch. Ration. Mech. Anal.}, {\bf 82}(1983), 313-346.

\bibitem{Berestycki-Lions=1983-2}
H. Berestycki and P. L. Lions. 
Nonlinear scalar field equations II: Existence of infinitely many solutions. 
{\it Arch. Ration. Mech. Anal.}, {\bf 82}(1983), 347-375.


\bibitem{Jeanjean=TransAMS=2019}
D.~Bonheure, J.~Casteras, T.~Gou and L.~Jeanjean. 
Normalized solutions to the mixed dispersion nonlinear Schr\"odinger equation in the mass critical and supercritical regime.
{\it Trans. Amer. Math. Soc.}, {\bf 372}(2019), no.3, 2167–2212.

 

\bibitem{Cazenave}
T. Cazenave. 
Semilinear Schr\"odinger Equations. 
Courant Lecture Notes in Mathematics, vol.10, New York University, New York, 2003.



\bibitem{Chang=2005}
Kung-Ching, Chang. 
Methods in nonlinear analysis. 
Springer Monographs in Mathematics. Springer-Verlag, Berlin, 2005.


  
 
 \bibitem{Colin-Jeanjean=NA=2004}
M. Colin and L. Jeanjean. 
Solutions for quasilinear Schr\"odinger equation: a dual approach. 
{\it Nonlinear Anal.}, {\bf 56}(2004), 213–226.

 
\bibitem{Colin-Jeanjean-Squassina=2010}
M. Colin, L. Jeanjean and M. Squassina. 
Stability and instability results for standing waves of quasi-linear Schr\"odinger equations.
{\it Nonlinearity}, {\bf 23}(2010), no.6, 1353–1385.


\bibitem{Coron=1984}
J.~M.~Coron. 
The continuity of the rearrangement in $W^{1,p}(\RN)$. 
{\it Ann. Scuola Norm. Sup. Pisa Cl. Sci. (4)}, {\bf 11}(1984), no.1, 57-85.


  
 
 
\bibitem{Ghoussoub=1993}
N.~Ghoussoub. 
Duality and perturbation methods in critical point theory, {\it Research Monograph, Cambridge Tracts, Cambridge University Press,} (1993)  268pp. 
 

\bibitem{GouJeanjean=2018}
T.~Gou and L.~Jeanjean.
\newblock Multiple positive normalized solutions for nonlinear Schr\"odinger systems.
{\it Nonlinearity}, {\bf 31}(2018), no.5, 2319-2345.



\bibitem{GouJeanjean=2016}
T.~Gou and L.~Jeanjean.
\newblock Existence and orbital stability of standing waves for nonlinear Schr\"odinger systems.
{\it Nonlinear Anal.}, {\bf 144}(2016), 10-22.


\bibitem{bg-7}
W. Hasse. 
A general method for the solution of nonlinear soliton and kink Schr\"odinger equations. 
{\it Z. Phys. B}, {\bf 37}(1980), 83-87.

 


 
 \bibitem{Ikoma=AdvDE=2019}
N.~Ikoma, K.~Tanaka. 
A note on deformation argument for $L^2$ normalized solutions of nonlinear Schr\"odinger equations and systems. 
{\it Adv. Differential Equations}, {\bf 24}(2019), 609-646.


 
\bibitem{Jeanjean=1997}
L.~Jeanjean. 
Existence of solutions with prescribed norm for semilinear elliptic equations. 
{\it Nonlinear Anal.}, {\bf 28}(1997), no.10, 1633-1659.


  
 \bibitem{Jeanjean-Lu=CVPDE=2020}
L. Jeanjean and S. S. Lu. 
A mass supercritical problem revisited. 
{\it Calc. Var. Partial Differential Equations}, {\bf 59}(2020), no.5, 44 pp.


\bibitem{Jeanjean-Luo=2013}
L. Jeanjean and T. J. Luo. 
Sharp non-existence results of prescribed $L^2$-norm solutions for some class of Schr\"odinger-Poisson and quasi-linear equations.
{\it Z. Angew. Math. Phys.}, {\bf 64}(2013), 937–954.
   
\bibitem{Jeanjean-Luo-Wang=JDE=2015}
L. Jeanjea, T. Luo and Z. Q. Wang. 
Multiple normalized solutions for quasi-linear Schr\"odinger equations. 
{\it J. Differential Equations}, {\bf 259}(2015), 3894-3928.

 

\bibitem{bg-2}
M. Kosevich, A. Ivanov and S. Kovalev. 
Magnetic solitons. 
{\it Phys. Rep.}, {\bf 194}(1990), 117-238.

 
 \bibitem{bg-4}
S. Kurihara. 
Large-amplitude quasi-solitons in superfluid films. 
{\it J. Phys. Soc. Jpn.}, {\bf 50}(1981), 3262–3267.



\bibitem{Kuzin-Pohozaev=1997}
I. Kuzin and S. Pohozaev.  
Entire solutions of semilinear elliptic equations. 
Progress in Nonlinear Differential Equations and their Applications, 33. Birkhäuser Verlag, Basel, 1997.



  
\bibitem{Li-Zou}
H. W. Li and W. M. Zou. 
Normalized ground states for semilinear elliptic systems with critical and subcritical nonlinearities. 
https://arxiv.org/abs/2006.14387.



\bibitem{MathN=2019}
Q. Q. Li, W. B. Wang, K. M. Teng and X. Wu. 
Multiple solutions for a class of quasilinear Schr\"odinger equations. 
{\it Math. Nachr.}, {\bf 292}(2019), no.7, 1530-1550.


 
 \bibitem{Loss-Lieb}
E. H. Lieb and M. Loss. 
Analysis. Second edition. 2001.



\bibitem{Lions2}
P.~L.~Lions. 
The concentration-compactness principle in the calculus of variations. The locally compact case. II.
{\it Ann. Inst. H. Poincaré Anal. Non Linéaire}, {\bf 1}(1984), no.4, 223–283.
 


\bibitem{bg-5}
G. Litvak and M. Sergeev. 
One dimensional collapse of plasma waves. 
{\it JETP Lett.}, {\bf 27}(1978), 517–520.



\bibitem{Liu-Liu-Wang=CVPDE=2013}
X. Q. Liu, J. Q. Liu and Z. Q. Wang. 
Ground states for quasilinear Schr\"odinger equations with critical growth.
{\it Calc. Var. Partial Differential Equations}, {\bf 46}(2013), 641-669.

 
\bibitem{Liu-Liu-Wang=ProAMS=2013}
X. Q. Liu, J. Q. Liu and Z. Q. Wang.  
Quasilinear elliptic equations via perturbation method. 
{\it Proc. Amer. Math. Soc}, {\bf 141}(2013), 253–263.


\bibitem{Liu-Liu-Wang=JDE=2013}
X. Q. Liu, J. Q. Liu, and Z. Q. Wang.  
Quasilinear elliptic equations with critical growth via perturbation method. 
{\it J. Differential Equations}, {\bf 254}(2013), 102–124.

\bibitem{Liu-Wang-Wang=JDE=2003}
J. Q. Liu, Y. Q. Wang and Z. Q. Wang. 
Soliton solutions for quasilinear Schr\"odinger equations II.  
{\it J. Differential Equations}, {\bf 187}(2003), 473–493.

\bibitem{Liu-Wang=JDE=2014}
J. Q. Liu and Z. Q. Wang.  
Multiple solutions for quasilinear elliptic equations with a finite potential well. 
{\it J. Differential Equations}, {\bf 257}(2014), 2874–2899.
 


\bibitem{bg-8}
G. Makhankov and K. Fedyanin. 
Non-linear effects in quasi-one-dimensinal models of condensed matter theory. 
{\it Phys. Rep.}, {\bf 104}(1984), 1-86.



\bibitem{GNinequality}
L. Nirenberg.   
On elliptic partial differential equations.  
{\it Ann. di Pisa}, {\bf 9}(1962), 187-195.


 
 \bibitem{Palais=1979}
R. S. Palais.  
The principle of symmetric criticality.  
{\it Commun. Math. Phys.}, {\bf 69}(1979), 19-30.

 
 

\bibitem{Poppenberg-Schmitt-Wang=CVPDE=2002}
M. Poppenberg, K. Schmitt and Z. Q. Wang.  
On the existence of soliton solutions to quasilinear Schr\"odinger equations. 
{\it Calc. Var. Partial Differential Equations}, {\bf 14}(2002), 329–344.


\bibitem{bg-6}
M. Porkolab and V. Goldman. 
Upper hybrid solitons and oscillating two-stream instabilities. 
{\it Phys. Fluids.}, {\bf 19}(1976), 872–881.




\bibitem{bg-3}
W. Quispel and W. Capel. 
Equation of motion for the Heisenberg spin chain. 
{\it Phys. A.}, {\bf 110}(1982), 41-80.


 \bibitem{Rabinowitz=1986}
P. H. Rabinowitz. 
Minimax Methods in Critical Point Theory with Applications to Differential Equations. 
CBMS Regional Conference Series in Mathematics, Vol. 65. American Mathematical Society, Providence, 1986.
 

\bibitem{Serrin-Tang=2000}
J. Serrin and M. Tang. 
Uniqueness of ground states for quasilinear elliptic equations.
{\it Indiana Univ. Math. J. }, {\bf 49}(2000), no.3, 897-923.


\bibitem{Soave=JDE=2020}
N.~Soave.
\newblock Normalized ground states for the NLS equation with combined nonlinearities.
{\it J. Differential Equations}, {\bf 269}(2020), no.9, 6941-6987.

\bibitem{Soave=JFA=2020}
N.~Soave.
\newblock Normalized ground states for the NLS equation with combined nonlinearities: the Sobolev critical case.
{\it J. Funct. Anal.}, {\bf 279}(2020), 108610.

\bibitem{Szulkin-Weth=JFA=2009}
A. Szulkin and T. Weth. 
Ground state solutions for some indefinite variational problems. 
{\it J. Funct. Anal.}, {\bf 55}(2009), 3802-3822.



 \bibitem{Ye=JMAA=2021}
H.~Y.~Ye and Y.~Y.~Yu. 
The existence of normalized solutions for $L^2$-critical quasilinear Schr\"odinger equations. 
{\it J. Math. Anal. Appl.}, {\bf 497}(2021), no.1, 124839.



\bibitem{Zeng-Zhang=2018}
X. Y. Zeng and Y. M. Zhang. 
Existence and asymptotic behavior for the ground state of quasilinear elliptic equations. 
{\it Adv. Nonlinear Stud.}, {\bf 18}(2018), no.4, 725–744.
 


\end{thebibliography}
 \end{document}